\documentclass[11pt]{article}
 
\usepackage{latexsym}
\usepackage{amsmath}
\usepackage{amsfonts}
\usepackage{amssymb}
\usepackage{graphpap}
\usepackage{epsfig}
\usepackage{amscd}
\usepackage{amsthm}
\usepackage{eucal} 

\font\got=eufm10 at 12pt
\font\posebni=msam10

\newcommand{\nor}[1]{|\hskip -0.6pt | #1 |\hskip -0.6pt |}

\newcommand{\mn}[2]{\{ #1\, ;#2 \}}
\newcommand{\sk}[2]{\langle #1 , #2\rangle}
\newcommand{\av}[1]{\langle #1 \rangle}

\newcommand{\dok}{\noindent {\bf Proof.}\ }
\newcommand{\Do}{\noindent {\bf Proof.}\ }

\newcommand{\C}[0]{{\mathbb C}}
\newcommand{\N}[0]{{\mathbb N}}

\newcommand{\R}[0]{\mathbb{R}}
\newcommand{\RR}[0]{{\mathbb R}^2}
\newcommand{\Z}[0]{{\mathbb Z}}
\newcommand{\PP}[0]{\mathbb{P}}
\newcommand{\EP}[0]{\mathbb{EP}}

\newcommand{\LL}[0]{{\cal L}}

\newtheorem{thm}{Theorem}
\newtheorem{lema}{Lemma}
\newtheorem{prop}{Proposition}
\newtheorem{cor}{Corollary}

\title{
Sharp estimates of martingale transforms in higher dimensions and applications to the Ahlfors-Beurling operator}
\author{\sc Oliver Dragi\v{c}evi\'c\footnote{ 
Supported by the European Comission (IHP network Harmonic Analysis and Related Problems 2002-2006, contract no. HPRN-CT-2001-00273-HARP) 
and the Ministry of Higher Education, Science and Technology of Slovenia (research program Analysis and Geometry, contract no. P1-0291).}, \ \sc Stefanie Petermichl \\
\sc and Alexander Volberg}
\date{March 2, 2006}

\begin{document}

\maketitle


\section{Introduction}

The main aspiration of this note is to construct several different Haar-type systems in euclidean spaces of higher dimensions and prove sharp $L^p$ bounds for the corresponding martingale transforms. In dimension one this was a result of Burkholder and in our proof we also rely on his lemma for abstract martingales involving differential subordination.

The motivation for working in this direction is the search for $L^p$ estimates of the {\it Ahlfors-Beurling operator} $T$, which is defined by
$$
Tf(z)=-\frac{1}{\pi}  \,{\rm p.v.}
\int_{\C}\frac{f(\zeta)}{(z-\zeta)^2}\ dA(\zeta)\,.
$$
There is a long-standing conjecture by T. Iwaniec \cite{Iw1} which asserts that 
\begin{equation}
\label{domneva}
\nor{T}_{B(L^p(\C))}=p^*-1\,,
\end{equation}
where $p^*=\max \{p,q\}$ and $q$ is the conjugate exponent of $p$. Note that
$$
p^*-1=\left\{
\begin{array}{lrr}
\ p-1 & ; & p\geqslant 2\\
(p-1)^{-1} & ; & 1<p\leqslant 2
\end{array}
\right.\,.
$$
While the conjecture is yet unconfirmed, it is known that the growth of norms is indeed linear. 
In \cite{DV1} the estimate 
$$
\nor{T}_p\leqslant\sqrt{2}(p-1)\Big(\frac1{2\pi}\int_0^{2\pi}|\cos\vartheta|^p\,d\vartheta\Big)^{-1/p}\,,\hskip 10pt p\geqslant 2\,,
$$
was given. Very recently a better estimate $\nor{T}_p\leqslant\sqrt{2p(p-1)}$ was obtained in \cite{BJ}, which under interpolation improves to $1.575(p^*-1)$. For large $p$ both \cite{DV1} and \cite{BJ} return $\text{\posebni\char47}\ \sqrt2(p-1)$.

In \cite{DV} the Ahlfors-Beurling operator $T$ was represented, up to some absolute constant $C$, as an average of two-dimensional martingale transforms. Our intention was to generalize the Burkholder's sharp estimate \cite{Bu1} of martingale transforms on $\R$ and apply it to this in order to obtain a new estimate on $T$. Although at the end we do not improve the best known constant for $\nor{T}_p$, the method applied turns out to be very useful for estimates of the powers of $T$ (see \cite{DPV}). Also, we present complete and explicit calculations of formulas for constants $C$ in each of these new systems. Sometimes, though, we were not able to give numerical value to these constants, so we do not know whether they would contribute any new estimates of $\nor{T}_p$. 

Another novelty of this work is that in order to obtain the operator $T$, we need to consider only one (very simple) martingale transform. See discussion on p. \pageref{concrete}. We believe that if one could estimate this concrete martingale transform by a constant smaller that $p^*-1$, this might reflect in a new $L^p$ bound for the Ahlfors-Beurling operator.

\bigskip

The theory giving rise to the conjecture \eqref{domneva} is extensive and comprises work of several generations of mathematicians. The sheer boundedness of $T$ on $L^p(\C)$, $1<p<\infty$, is a consequence of the Calder\'on-Zygmund theory. What distinguishes $T$ from other similar operators is that it interchanges $\partial$ and $\bar\partial$ derivatives: $T(\bar\partial f)=\partial f$ for $f$ in the Sobolev class $W^{1,2}(\C)$. This is the source of its significance for nonlinear PDE and quasiconformal mappings on $\C$. 
A $K$-$quasiconformal$ $mapping$ can be defined as a generalized $L^2$-solution of the $Beltrami$ $equation$ $\bar\partial f=\mu\cdot\partial f$, where $\nor{\mu}_\infty=
\frac{K+1}{K-1}<1$.
It is known since the early work of Bojarski \cite{Bo}, \cite{Bo1} that quasiconformal maps belonging to $W^{1,2}_{loc}$ in fact also belong to $W^{1,p}_{loc}$ for some $p>2$ depending only on $\nor{\mu}_\infty$. This property relates closely to area distortion of quasiconformal mappings. A great deal of effort was devoted to establishing the best integrability, i.e. the limits of the said self-improvement. This problem was settled in 1994 by Astala \cite{A}. A ``dual" result was confirmed by Petermichl and Volberg in \cite{PV}. 
The Ahlfors-Beurling operator plays an important part here, e.g. see \cite{LV}. 
Also the conjecture \eqref{domneva} was based on circumstancial evidence arising from this theory \cite{Iw1}. In particular, if true, it would imply Astala's area distortion theorem. For a wider discussion we point, for example, to \cite{BM-S}, where it is also touched upon the relation between $T$ and the Morrey's problem concerning quasiconvexity.
 
\bigskip
From now on, we will use the symbol $\nor{\cdot}_p$ to denote both $L^p$ norms of functions as well as operator norms of operators acting on $L^p$, unless specified otherwise. The term {\it test function} will stand for a continuous, compactly supported, complex function on $\R$ (or $\R^n$, more generally) whose average is zero: $\int_{\R^n}f(x)\,dx=0$. Such functions are dense in $L^p(\R^n)$ for $1<p<\infty$.

\subsection{Main idea}
Our approach to this problem is based on the principal result of \cite{DV}. There it was proven that $T=C\cdot T'$, where $T'$ was a result of an averaging process of two-dimensional martingale transforms. This equality was then applied 
on  weighted $L^p$ spaces, 
but this time we were aiming at the unweighted case. Namely, Burkholder \cite{Bu2} showed that every one-dimensional martingale transform $T_\sigma$ admits an estimate $\nor{T_\sigma}_p\leqslant p^*-1$, which is sharp. 
Moreover, a closer examination shows that 
$$
\aligned
\label{izraz}
C & = \frac{\ \hskip 12pt 12\log2}{16\pi + 32\log 2 - 15 \log 5 - 40\arctan 2}\\
& \ \text{\posebni\char47}\
2,07\,.
\endaligned
$$
Thus the idea is to 
\begin{itemize}
\item
extend the Burkholder's estimate in order to obtain $p^*-1$ estimates for {\sl two}--dimensional martingale transforms 
\item
reduce them to the estimate of one simple martingale transform which, if tamed, could give an improvement of the known $L^p$ estimates for $T$.
\end{itemize}

\subsection{Statement of results}

We begin by revisiting Burkholder's theorem in Section \ref{mt}. 
After treating the linear case, we define a two-dimensional Haar base (and call it ${\cal H}_{orig}$) the way it was done in \cite{DV} and attempt applying the same proof to this case, as well. There we explain why generalizing the Burkholder's theorem is not straightforward. The encountered obstacle is removed 
in the continuation of the section by means of constructing ${\cal H}_{new}$, the first of our modified Haar systems. It is derived from ${\cal H}_{orig}$ by introducing an apparently small, but vital modification. 

In ${\cal H}_{new}$ we also run a process of averaging martingale transforms which at the end gives a constant multiple of $T$. Together with the obtained generalization of the Burkholder's theorem this means a $C(p^*-1)$ bound for $T$. The size of $C$ is approximately 2,07 and it appeared already after the averaging process in \cite{DV}. In Sections \ref{ponoc} and \ref{largo} we concentrate on improving it. For that purpose we gradually define several further classes of Haar systems on the plane, called ${\cal H}_{b,\varphi}$, ${\cal H}_{b,\varphi}^d$ and ${\cal H}_{a,b}^\vartriangle$, all of which carry on certain geometric structure.

Subsection \ref{hidi} brings a short account on how our estimates could be extended to arbitrary $\R^n$ where they might work for some homogeneous kernels of degree $-n$.
By ${\cal H}_{new}^n$ we will denote a generalization of ${\cal H}_{new}$ in $\R^n$. 
We summarize our results as follows.

\bigskip
\noindent
{\bf Theorem} 
{\it 
Let ${\cal H}$ be any of the Haar systems 
${\cal H}_{b,\varphi}$, ${\cal H}_{b,\varphi}^d$, ${\cal H}_{a,b}^\vartriangle$ (in $\R^2$) or ${\cal H}_{new}^n$ (in $\R^n$). If $T_\sigma$ is any martingale transform associated to ${\cal H}$, then $\nor{T_\sigma}_p\leqslant p^*-1$. This estimate is sharp.}  

\bigskip
\noindent
{\bf Theorem} 
{\it Let ${\cal H}$ be any of the Haar systems 
${\cal H}_{b,\varphi}$, ${\cal H}_{b,\varphi}^d$ and ${\cal H}_{a,b}^\vartriangle$. Then $T=C\cdot T'$, where $T'$ is obtained as an average of martingale transforms in ${\cal H}$. The best $C$ is not larger than approximately $2,007$. 
Consequently, $\nor{T}_p\ \text{\posebni\char47}\ 2,007(p^*-1)$}. \\ 

In Sections \ref{ponoc} and \ref{largo} we present detailed calculations of constants appearing in the Theorem above.




\section{
Sharp estimates for martingale transforms}
\label{mt}

In order to extend the Burkholder's result it is useful to understand how it was obtained in the first place. A short, elegant proof is described in \cite{Bu2}. However, there he gives no explanation on how he obtained his {\it Bellman function}, which lies at the heart of all known estimates of the $p^*-1$ type. This gap is filled in his exposition \cite{Bu}. In that work, on p. 16, he also proves the following lemma, which is very useful due to its generality and sharpness. We present it for the convenience of the reader. For the details of the proof one should refer to \cite{Bu3}.

\begin{lema}
\label{Bu}
Let $(${\got W}$, {\cal F}, P)$ be a probability space, $\mn{{\cal
F}_n}{n\in\N}$ a filtration in ${\cal F}$ and $H$ a separable Hilbert space. Furthermore, let
$(X_n, {\cal F}_n, P)$ and $(Y_n, {\cal F}_n, P)$ be $H$-valued
martingales satisfying
\begin{equation}
\label{difsub}
\|Y_0(\omega)\|_H \leq \|X_0(\omega)\|_H\ \text{and }
\|Y_n(\omega)-Y_{n-1}(\omega)\|_H \leqslant
\|X_n(\omega)-X_{n-1}(\omega)\|_H
\end{equation}
for all $n\in\N$ and almost every $\omega\in$ {\got W}. Then for
any $p\in (1,\infty)$
$$
\|Y_n\|_p
\leqslant (p^*-1) \|X_n\|_p
\,.
$$
The constant $p^*-1$ is sharp.
\end{lema}

The property \eqref{difsub} is called {\it differential subordination}. 
%
In order to proceed to the proof of the sharp estimate for one-dimensional martingale transforms, we require some further notation.

\bigskip

We call the family of intervals ${\cal L}:=\mn{[m\, 2^n, (m+1)\,
2^n)}{m,n\in\Z}$ the standard {\it dyadic lattice}. 
Each interval $I\subset \R $ gives rise to its
{\it Haar function} $h_I$, defined by
$$
h_I:=|I|^{-1/2}(\chi_{I_+}-\chi_{I_-}) \, ,
$$
where $I_-$ and $I_+$ denote the left and the right half of
the interval $I$ respectively, and $\chi_E$ stands for the
characteristic function of the set $E$, as usual. 
Denote by ${\cal L}(I)$ the set of all dyadic subintervals of the interval $I$, including $I$ itself.
For any $p\in(1,\infty)$ and any interval $I$, the set $\mn{h_J}{J\in{\cal L}(I)}$ 
forms a {\it basis} of the space $L^p(I)$. By that we shall mean that for $f\in L^p(I)$,
$$
f-\av{f}_I\chi_I=\lim_{n\rightarrow\infty}\sum_{J\in {\cal L}(I)\atop |J|>2^{-n}|I|}\sk{f}{h_J}h_J\,,
$$
the limit existing in the $L^p$-sense and $\av{f}_I$ standing for the average of function $f$ over $I$. A similar statement is valid for arbitrary intervals, of course.

Now we are able to define
the operator $T_{\sigma}$ by
$$
T_{\sigma}f:=\sum_{J\in {\cal L}}
\sigma_J\sk{f}{h_J}\,
h_J\, ,
$$
where $\sigma : {\cal L} \rightarrow S^1$ is
arbitrary. Such operators are called {\it martingale transforms}. Note that if $f$ is a test function, the terms $\sk{f}{h_J}$ are nonzero only for $J$ contained in the support of $f$.

\begin{thm}
\label{complex-functions}
For every $\sigma$ and every $1<p<\infty$,
$$
\nor{T_{\sigma}}_{p} \leqslant p^*-1\,.
$$
\end{thm}

\dok
The present argument is somewhat different from the original one, published in \cite{Bu2}, 
in that it applies the above lemma rather than the Bellman function of Burkholder directly. 

It suffices to start with a test function $f$ with its support contained in  $I\in\mathcal{L}$. Recall that this implies $\av{f}_I=0$.
Take $H=\C$ and let ${\cal F}_n$ be the $\sigma$-algebra of all
dyadic subintervals of $I$ with length at least $2^{-n}|I|$. For $n\in\N\cup\{0\}$ and 
$\omega\in I$ define
$$
X_n (\omega):= 
\sum_{J\in {\cal L}(I) 
\atop |J|>2^{-n}|I|}
\sk{f}{h_J}h_J(\omega)\,.
$$
Choose a sequence of numbers $\sigma_J \in S^1$
and consider
$$
Y_n (\omega):= \sum_{J\in {\cal L}(I) 
\atop |J|> 2^{-n}|I|}
\sigma_J\sk{f}{h_J}h_J(\omega)\,.
$$
It follows immediately from the construction that both $(X_n, {\cal F}_n,
dx)$ and $(Y_n, {\cal F}_n, dx)$ are martingales. Clearly
$$
X_{n+1}-X_{n} = \sum_{J\in {\cal L}(I) 
\atop |J|= 2^{-n}|I|} \sk{f}{h_J}h_J
$$
as well as
$$
Y_{n+1}-Y_{n}= \sum_{J\in{\cal L}(I)
\atop |J|= 2^{-n}|I|}
\sigma_J\sk{f}{h_J}h_J\,.
$$
Since $|\sigma_J|=1$ and the sums above comprise functions with disjoint supports,
\hyphenation{interiors}
these martingales are differentially subordinated to each other, so we may apply Lemma \ref{Bu}. Since $\|f\|_{p} =
\lim_{n\rightarrow\infty}\|X_n\|_{p}$ and $\|T_{\sigma}f\|_{p}
= \lim_{n\rightarrow\infty}\|Y_n\|_{p}$, 
we completed the proof.
\qed

\subsection{Two dimensional case}

We would like to extend this result to the martingale transforms on the plane. As mentioned before, we want to do this by applying again Lemma \ref{Bu}. 

Let us start with the construction of the Haar system on the plane. We repeat the definitions from \cite{DV}, since the main theorem there (Theorem 1), where $T$ is  represented as an average of planar martingale transforms, is the one we are aiming to utilize.
\label{plane}

\bigskip

The term {\it dyadic lattice} and the symbol ${\cal L}$ will now stand for
the collection of all squares of the form $I\times J \subset
\R^2$, where $I$ and
$J$ are dyadic intervals of the same length. To each such
square $Q=I\times J$ we will assign three Haar functions:
$$
\aligned
& h_Q^1(s,t)=\chi_I(s)h_J(t)|I|^{-1/2}\\
& h_Q^2(s,t)=h_I(s)\chi_J(t)|J|^{-1/2}\\
& h_Q^3(s,t)=h_I(s)h_J(t)\,.
\endaligned
$$
\noindent
Symbolically,
%
%
%
%
$$ 
h_Q^1 \equiv
\begin{array}{|c|}
\hline \raisebox{0pt}[12pt][6pt]{\hskip 10pt$+$\hskip 10pt}
\\
\hline \raisebox{0pt}[12pt][6pt]{\hskip 10pt $-$\hskip 10pt}   \\
\hline
\end{array}
\hskip 30pt 
h_Q^2 \equiv
\begin{array}{|c|l|}
\hline \raisebox{0pt}[21pt][15pt]{$-$} & + \\ \hline
\end{array}
\hskip 30pt 
h_Q^3 \equiv
\begin{array}{|c|c|}
\hline \raisebox{0pt}[12pt][6pt]{$-$} & +  \\ \hline
\raisebox{0pt}[12pt][6pt]{$+$} &  -  \\ \hline
\end{array}
$$
As previously, one can verify that the set $\mn{h_Q^i}{Q\in {\cal L},\ i=1,2,3}$ constitutes a  basis of $L^p(\R^2)$. In order to distinguish it from the subsequent Haar systems, we will call it ${\cal H}_{orig}$. 
Now the {\it two-dimensional martingale
transform} becomes the operator
$$
T_{\sigma}f:=\sum_{Q\in {\cal L}}\sum_{i=1}^3
\sigma^i_Q\sk{f}{h_Q^i}\,h_Q^i \, ,
$$
where, as before, $\sigma^i : {\cal L}
\rightarrow S^1$. 

\subsection{The problem}
\label{theproblem}

Let us try to repeat the proof that worked in the linear case. Take a test function $f$ be supported on some dyadic square $\Omega$. 
We may assume for simplicity that $\Omega$ is of size 1.
Similarly as before, let ${\cal L}(\Omega)$ stand for the set of all dyadic subsquares of $\Omega$. Denote, for $n\in\N$,
$$ 
X_n=
\sum_{Q\in {\cal L}(\Omega)\atop |Q|> 4^{-n}}
\sum_{i=1}^3
\sk{f}{h_Q^i}\,h_Q^i 
$$
and
$$ 
Y_n=\sum_{Q\in {\cal L}(\Omega)\atop |Q|> 4^{-n}}\sum_{i=1}^3
\sigma_Q^i\sk{f}{h_Q^i}\,h_Q^i \, .
$$
If ${\cal F}_n$ is the $\sigma-$algebra of squares
in ${\cal L}(\Omega)$ of the size at least $4^{-n}$, then $(X_n,{\cal
F}_n)$ and $(Y_n,{\cal F}_n)$ are martingales with respect to the
planar Lebesgue measure.
But now the (martingale) differences
$$
X_{n+1}-X_n=\sum_{Q\in {\cal
L}(\Omega)\atop |Q|= 4^{-n}} 
\sum_{i=1}^3
\sk{f}{h_Q^i}\,h_Q^i \, 
$$
and
$$
Y_{n+1}-Y_n =\sum_{Q\in {\cal L}(\Omega)\atop |Q|= 4^{-n}}
\sum_{i=1}^3
\sigma_Q^i\sk{f}{h_Q^i}\,h_Q^i \, $$

\noindent
are not sums of functions with mutually disjoint
supports, like in the one-dimensional case, because this time we have
three Haar functions associated to each square, and not just one,
which quickly implies that the differential subordination from Burkholder's 
Lemma \ref{Bu} does {\sl not} hold anymore, unless in the trivial case.

In order to fix this problem, when passing from $X_m$ to $X_{m+1}$ we could take just one type of Haar functions associated to the smaller squares, rather than three. For example, if, as before, 
$$ 
X_{3n}=\sum_{Q\in {\cal L}(\Omega)\atop |Q|> 4^{-n}}
\sum_{i=1}^3
\sk{f}{h_Q^i}\,h_Q^i \, ,
$$
then $X_{3n+1}$ could be defined as 
$$ 
X_{3n+1}=X_{3n}+
\sum_{Q\in {\cal L}(\Omega)\atop |Q|=4^{-n}}
\sk{f}{h_Q^1}\,h_Q^1 
$$
and similarly
$$ 
X_{3n+2}=X_{3n+1}+
\sum_{Q\in {\cal L}(\Omega)\atop |Q|=4^{-n}}
\sk{f}{h_Q^2}\,h_Q^2 \, .
$$
If we define $Y_m$ analogously, then we managed to remove the previous obstacle, that is, we obtain
$$
|Y_{m+1}-Y_m|\leqslant |X_{m+1}-X_m|
$$
pointwise on $\C$ (in fact, we obviously have equality). However, we acquired another problem. 

For these modified sequences $X$ and $Y$ to be martingales, we have to specify the filtration $\mn{{\cal F}_m}{m\in\N}$. For ${\cal F}_{3n}$ we take the sub-$\sigma$-algebra generated by all dyadic squares in $\Omega$ of size $4^{-n}$. Following the nature of $h_Q^1$ we let ${\cal F}_{3n+1}$ be generated by the upper and lower halves of squares from ${\cal F}_{3n}$. Similarly, ${\cal F}_{3n+2}$ is generated by ${\cal F}_{3n+1}$ and left and right halves of squares in ${\cal F}_{3n}$. But the problem is that now ${\cal F}_{3n+2}$ is equal to ${\cal F}_{3n+3}$. Consequently, the conditional expectation $\mathbb{E}(X_{3n+3}|{\cal F}_{3n+2})$ is equal to $X_{3n+3}$ and not $X_{3n+2}$. In other words, $(X_m,{\cal F}_{m},dx)$ and $(Y_m,{\cal F}_{m},dx)$ are two stochastic processes which mutually satisfy differential subordination, but they are not martingales, so we again cannot apply the Burkholder's lemma.
\label{problem}

Note that the described obstacle occurs (only) on every third step. Hence when passing from one sub-$\sigma$-algebra, generated by full squares, to another (in our notation, from ${\cal F}_{3n}$ to ${\cal F}_{3n+3}$), we can afford only one intermediate step, but instead we have two, thus we should somehow eliminate one. 

One way is to do that ``artificially". That is, if we take functions which are orthogonal to all $h_Q^i$, $Q\in{\cal L}(\Omega)$, for some (fixed) $i\in\{1,2,3\}$, we get two martingales and can apply Lemma \ref{Bu}. Thus the restriction of $T_\sigma$ to any of the 
subspaces $\mn{h_Q^i}{Q\in{\cal L}(\Omega)}^\bot$, $i=1,2,3$, is of $L^p$ norm which does not exceed $p^*-1$. 


\subsection{Remedy}
\label{remedy}

There is a way out of the problems encountered. It was suggested to us by Guy David. 
His solution is simple yet very effective. It lies in changing the (Haar) system, for there is no reason to consider ${\cal H}_{orig}$ privileged over other possible bases. Let us explain how to do that.

To each square $Q$ we associate a different set of Haar functions:
%
\begin{equation}
\label{new}
\aligned
& h_{Q_0}:=h_Q^1\\
& h_{Q_+}:=\frac{1}{\sqrt{2}}(h_Q^2+h_Q^3)\\
& h_{Q_-}:=\frac{1}{\sqrt{2}}(h_Q^2-h_Q^3)\,.
\endaligned
\end{equation}
\noindent
Symbolically,
%
%
%
%
%
$$ 
h_{Q_0} \equiv
\begin{array}{|c|}
\hline \raisebox{0pt}[12pt][6pt]{\hskip 10pt$+$\hskip 10pt}
\\
\hline \raisebox{0pt}[12pt][6pt]{\hskip 10pt $-$\hskip 10pt}   \\
\hline
\end{array}
\hskip 30pt 
h_{Q_+} \equiv
\begin{array}{|c|c|}
\hline \raisebox{0pt}[12pt][6pt]{$-$} & +  \\ \hline
\multicolumn{2}{|c|}
{\raisebox{0pt}[12pt][6pt]{} }\\ \hline
\end{array}
\hskip 30pt 
h_{Q_-} \equiv
\begin{array}{|c|c|}
\hline 
\multicolumn{2}{|c|}
{\raisebox{0pt}[12pt][6pt]{} }   \\ \hline
\raisebox{0pt}[12pt][6pt]{$-$} &  +  \\ \hline
\end{array}
$$
Let us denote the system $\mn{h_{Q_*}}{*\in\{0,+,-\},\ Q\in\LL}$
by ${\cal H}_{new}$. 
Its big advantage is that $h_{Q_+}$ and $h_{Q_-}$ have disjoint supports. For in that case the associated martingale transforms do admit the desired estimates, as the following theorem shows.

\begin{thm}
\label{}
For any $Q\in \LL$ and $*\in\{0,+,-\}$ let $\sigma_{Q_*}$ be arbitrary unilateral complex numbers. Define the operator 
$$
T_\sigma f:=
\sum_{Q\in \LL}
\left[\sigma_{Q_0}\sk{f}{h_{Q_0}}h_{Q_0}
+\sigma_{Q_+}\sk{f}{h_{Q_+}}h_{Q_+}
+\sigma_{Q_-}\sk{f}{h_{Q_-}}h_{Q_-}\right].
$$
Then $\nor{T_\sigma}_{p}\leqslant p^*-1$. This estimate is sharp. 
\end{thm}

\dok
Take a test function $f$, supported in some $\Omega\in{\cal L}$, and 
define 
$$
\aligned 
X_{2n}& :=
\sum_{Q\in{\cal L}(\Omega)\atop |Q|>4^{-n}}
\left
[\sk{f}{h_{Q_0}}h_{Q_0}
+\sk{f}{h_{Q_+}}h_{Q_+}
+\sk{f}{h_{Q_-}}h_{Q_-}
\right]\\
X_{2n+1}& :=X_{2n}+
\sum_{Q\in{\cal L}(\Omega)\atop  
|Q|=4^{-n}}
\sk{f}{h_{Q_0}}h_{Q_0}
\endaligned
$$
and
$$
\aligned 
Y_{2n}& :=
\sum_{Q\in{\cal L}(\Omega)\atop  
|Q|>4^{-n}}
\left[\sigma_{Q_0}\sk{f}{h_{Q_0}}h_{Q_0}
+\sigma_{Q_+}\sk{f}{h_{Q_+}}h_{Q_+}
+\sigma_{Q_-}\sk{f}{h_{Q_-}}h_{Q_-}\right]\\
Y_{2n+1}& :=Y_{2n}+
\sum_{Q\in{\cal L}(\Omega)\atop 
|Q|=4^{-n}}
\sigma_{Q_0}\sk{f}{h_{Q_0}}h_{Q_0}\,.
\endaligned
$$
Let ${\cal F}_m$ be the $\sigma-$algebra, generated by $X_m$. Explicitly, ${\cal F}_{2n}$ is generated by all dyadic squares of size $4^{-n}$, while ${\cal F}_{2n+1}$ is generated by their upper and lower halves. This time, as opposed to the case described on page \pageref{problem}, ${\cal F}_{m+1}$ is properly contained in ${\cal F}_m$, hence $(X_m,{\cal F}_m,dx)$ and $(Y_m,{\cal F}_m,dx)$ are martingales. Moreover, it is clear that they satisfy the differential subordination: 
$$
|(X_{m+1}-X_m)(\omega)|=|(Y_{m+1}-Y_m)(\omega)| \hskip 20pt {\rm p.p.}\ \omega\in\C\,.
$$
We can apply Lemma \ref{Bu} and get that $\nor{Y_m}_p\leqslant (p^*-1)\nor{X_m}_p$ for every $m\in\N$. Now use that $\lim_{m\rightarrow\infty}\nor{X_m}_p=\nor{f}_p$ and $\lim_{m\rightarrow\infty}\nor{Y_m}_p=\nor{T_\sigma f}_p$.
\hfill\qed

\bigskip
\noindent
{\bf Remarks}
\label{sharpness}

--The result of the preceding theorem is sharp in the sense that 
$p^*-1$ is the supremum of $p$-norms of all martingale transforms ${T}_\sigma$.
This is a simple consequence of the fact that the same is true in Theorem  \ref{complex-functions} (see \cite{Bu}).

--Note that the supports of $h_{Q_+}$ and $h_{Q_-}$ being disjoint enabled us to add these functions in $X_n$'s and $Y_n$'s {\sl simultaneously} without risking any damage. This is how we lost the superfluous {\sl step} (cf. the end of the previous subsection). 

--This theorem can also be proven by 
introducing appropriate modifications to the Burkholder's direct proof of Theorem \ref{complex-functions}, which can be found in \cite{Bu2}. Of course, that would
not be the case had we replaced ${\cal H}_{new}$ by ${\cal H}_{orig}$.

\section{The averaging}
\label{ponoc}

Now that we constructed a Haar system for whose martingale transforms we are able to prove the $p^*-1$ estimates, we would like to apply this result to $T$. 
This section will be devoted to the proof of the following statement. 
%
%
\begin{equation}
\begin{array}{l}
\label{grobo}
\text{\it The Ahlfors-Beurling operator can be realized as a constant times an}\\
\text{\it average of two-dimensional martingale transforms.}
\end{array}
\end{equation}

This theorem, in its original form, first appeared in \cite{DV}. There it was proven for 
${\cal H}_{orig}$.
We are going to see that it actually holds for many different types of underlying Haar systems. Then our strategy will be to choose such that gives the best (i.e. smallest) constant, though always keeping in mind to choose among those which yield $p^*-1$ bounds. In order to do that, it will be necessary to have a scrupulous look at the proof of the theorem. For the convenience of the reader, we summarize it here as it appeared in \cite{DV}.

\bigskip

Instead of a dyadic lattice let us for a moment consider a unit
{\it grid} ${\cal G}$ of squares. This is a family of
squares 
$I\times J$, where $I$ and $J$ are
dyadic intervals of unit length. Furthermore, for $t\in
\R^2$  define ${\cal G}_t:={\cal G}+t$, i.e. the grid of
unit squares such that one of them contains
point $t$ as one of its vertices.

Introduce
$$
\PP_t f := \sum_{Q\in {\cal G}_t} \sk{f}{h_{Q_0}}h_{Q_0}\, .
$$
For the sake of transparency of the proof we only took $h_{Q_0}$ in the operator above. The procedure does not change  significantly if $h_{Q_+}$ and $h_{Q_-}$ are included.

Notice that the family $\Omega:=\mn{{\cal G}_t}{t\in\R^2}$
of all unit grids naturally corresponds to the torus
$\R^2/\Z^2$, which is of course in one-to-one
correspondence with the square
$[0,1)^2$. Thus we are able to regard $\Omega$ as a
probability space where the probability measure equals 
the Lebesgue measure on $[0,1)^2$.

Now this leads to the ``mathematical expectation" of the
``random variable" $\PP$. This will again be an operator on
$L^2(\C)$, defined pointwise (for $f\in L^2(\C)$) as
$$
(\EP\, f)(x):=\int_\Omega\PP_tf(x)\, dm(t)\, .
$$
Since $\EP$ is a result of integrating over a
certain probability space, 
it makes sense to call this process the {\it
averaging}. The structure of this operator is revealed
in the following proposition.

\begin{prop}
\label{F}
Assuming the notation as above, the operator $\EP$
is a convolution operator with the kernel
$
F(x_1,x_2)=-\beta(x_1)\alpha(x_2)\,,
$
where
$$
\alpha=h_0*h_0\ \ \ and\ \ \ \beta=\chi_0*\chi_0.
$$
Here $\chi_0$ and $h_0$ stand (respectively) for the characteristic and Haar function of the interval $[-1/2,1/2)$.
Inserting $h_Q^2$ in $\PP_t$ instead of $h_Q^1=h_{Q_0}$ yields $-\alpha(x_1)\beta(x_2)$, while $h_Q^3$ would produce $\alpha(x_1)\alpha(x_2)$. 
\end{prop}

\Do 
We are going to demonstrate the first part only, for the second one is proven in the same way. Choose $t=(t_1,t_2)\in\RR$, $Q=I\times J\in {\cal
G}_t$. Then
$$
\aligned
\sk{f}{h_{Q_0}} & = \int_\R\!\int_\R f(s_1,s_2)h_{Q_0}(s_1,s_2)\,
ds_1\, ds_2 
\,.
\endaligned
$$
Thus for (fixed) $f\in L^2(\R^2)$ and $x=(x_1,x_2)\in\R^2$ we have
$$
\aligned 
(\PP_t f)(x) & = \sum_{Q\in {\cal G}_t}\sk{f}{h_{Q_0}}h_{Q_0}(x)=\\
& =\sum_{Q\in {\cal G}_t}\int_\R\!\int_\R f(s_1,s_2)
h_{Q_0}(s_1,s_2)
\,ds_1\, ds_2 \, \cdot h_{Q_0}(x)\\
& =\int_\R\!\int_\R f(s_1,s_2)
\Big ( \sum_{Q\in {\cal G}_t}
\chi_I(s_1)h_J(s_2)h_{Q_0}(x)
\Big )
\,ds_1\, ds_2\, .
\endaligned
$$
The expression under the summation sign in the last row is
nonzero for exactly one $Q\in {\cal G}_t$; namely such that
$h_{Q_0}(x)\not =0$. This means that
$x=(x_1,x_2)\in Q$, therefore $x_1\in I$ and $x_2\in J$.
Thus $h_{Q_0}(x)=\chi_I(x_1)h_J(x_2)$. We get
\begin{equation}
\label{a}
(\PP_t f)(x)
=\int_\R\!\int_\R f(s_1,s_2)
\chi_I(s_1)h_J(s_2)\chi_I(x_1)h_J(x_2)  
\,ds_1\, ds_2\, .
\end{equation}

Because ${\cal G}_t$ does not change if we increase or
decrease any component of $t$ by 1, we may assume that
$I=[t_1-1,t_1)$ and $J=[t_2-1,t_2)$. Denoting
$I_0=[-1/2,1/2)$, this assumption implies
$$
I=t_1-\frac{1}{2}+I_0 \ \ \ {\rm and}\ \ \
J=t_2-\frac{1}{2}+I_0\, .
$$
Now let $\chi_0$ and $h_0$ be as in the formulation of the
Proposition. 
Then
$$
\chi_I(z)=\chi_0(z+1/2-t_1)
$$
and 
$$
h_I(z)=h_0(z+1/2-t_1)
$$
for all $z\in\RR$. The analogue pair is valid also for $J$,
of course.


Together with (\ref{a}), the last two equalities imply
$$
\aligned
(\PP_t f)(x) &=\int_\R\!\int_\R f(s_1,s_2)\,
\chi_0(s_1+1/2-t_1)h_0(s_2+1/2-t_2)\\
& \hskip 86pt \chi_0(x_1+1/2-t_1)h_0(x_2+1/2-t_2)
\,ds_1\,ds_2\, .
\endaligned
$$
Recall that $x_1\in I=[t_1-1,t_1)$ and $x_2\in
J=[t_2-1,t_2)$. Hence $x_i<t_i\leqslant x_i+1$ for $i=1,2$.
Averaging means in our case integrating over all admissible
$t_i$. Therefore
$$
(\EP f)(x)=\int_{x_2}^{x_2+1}\int_{x_1}^{x_1+1}
(\PP_{(t_1,t_2)}f)(x)\,dt_1\, dt_2\, .
$$
By using the most recent expression for $(\PP_tf)(x)$ and
changing variables (to $t_i-1/2$) we get
$$
(\EP f)(x)
=\int_{x_2-1/2}^{x_2+1/2}\int_{x_1-1/2}^{x_1+1/2}\
\int_\R\!\int_\R f(s_1,s_2) [A]\, ds_1\,ds_2\ dt_1\,
dt_2\, , 
$$
where
$$
A=\chi_0(s_1-t_1)h_0(s_2-t_2) \chi_0(x_1-t_1)h_0(x_2-t_2)\,.
$$
Applying Fubini's theorem yields
\begin{equation}
\label{b}
(\EP
f)(x_1,x_2)=
\int_\R\!\int_\R f(s_1,s_2)\
\int_{x_2-1/2}^{x_2+1/2}\int_{x_1-1/2}^{x_1+1/2}[A]\,
dt_1\, dt_2\ \
ds_1\,ds_2  \, .
\end{equation}
This is how we obtained the candidate for the
convolution kernel $F$ of the operator $\EP$. Namely, line
(\ref{b}) gives us the relation
$$
F(x_1-s_1,x_2-s_2)=\int_{x_2-1/2}^{x_2+1/2}\int_{x_1-1/2}^{x_1+1/2}[A]\,
dt_1\, dt_2\, .
$$
We are justified in writing this, since a change of variables shows that the above expression indeed depends only on $x_1-s_1$ and $x_2-s_2$. Besides, such a property is to be expected from a kernel, resulting from a process of averaging over translations.

Finally, taking $s_1=s_2=0$ returns
$$
\aligned
F(x_1,x_2)& =
\int_{x_2-1/2}^{x_2+1/2}\int_{x_1-1/2}^{x_1+1/2}
\chi_0(-t_1)h_0(-t_2)\chi_0(x_1-t_1)h_0(x_2-t_2)dt_1\, dt_2\, .\\
& = -\int_{x_1-1/2}^{x_1+1/2} 
\chi_0(t_1)\chi_0(x_1-t_1)\, dt_1
\cdot
\int_{x_2-1/2}^{x_2+1/2}h_0(t_2)h_0(x_2-t_2)\, dt_2\\
& = -(\chi_0*\chi_0)(x_1)\, (h_0*h_0)(x_2)\,,
\endaligned
$$
as desired. \hfill $\square$

\bigskip

Graphs of functions $\alpha$ and $\beta$ are shown as
Figures \ref{fig1} and \ref{fig2}, respectively.


\setlength{\unitlength}{1mm}
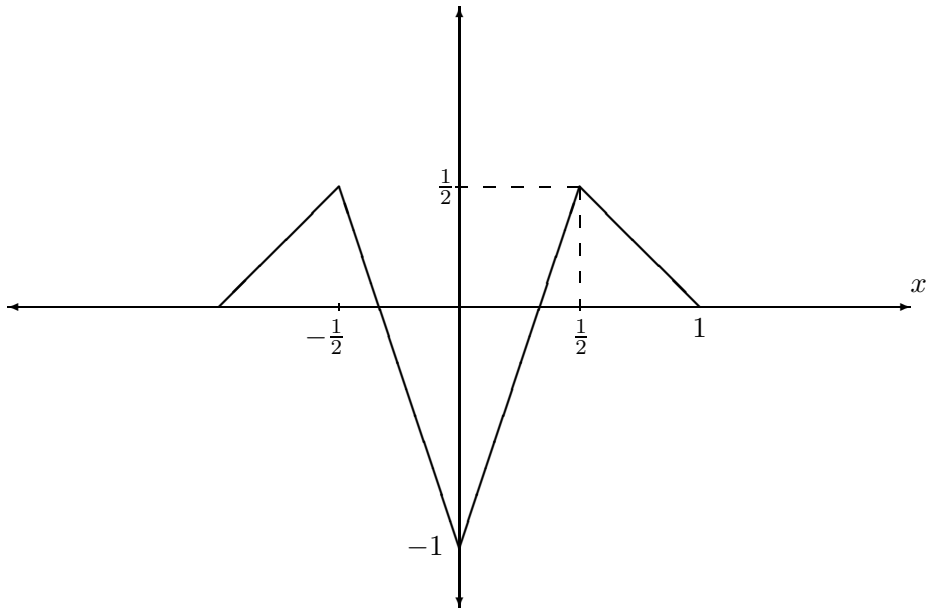
\begin{figure}
\begin{center}
\begin{picture}(140,80)(-68,-35)
\put(0,0){\vector(1,0){60}}
\put(0,0){\vector(-1,0){60}}
\put(0,0){\vector(0,1){40}}
\put(0,0){\vector(0,-1){40}}
\multiput(-0.5,16)(3.72,0){5}{\line(1,0){1.5}}
\multiput(16,-0.5)(0,3.7){5}{\line(0,1){1.5}}
\multiput(16,-0.5)(0,3.7){1}{\line(0,1){1}}
%
%
%
%
\put(15,-5){$\frac12$}
\put(-3,15){$\frac12$}
\multiput(-16,-0.5)(0,3.7){1}{\line(0,1){1}}
\put(-20.5,-5){$-\frac12$}

\put(60,2){$x$}
\put(31,-4){$1$}
\put(-7,-33){$-1$}
\put(-12,16){\line(1,-1){2}}

\thicklines
\put(0,-32){\line(1,3){16}}
\put(16,16){\line(1,-1){16}}
\put(0,-32){\line(-1,3){16}}
\put(-16,16){
\line(-1,-1){16}}


\end{picture}
\end{center}
\caption{Graph of $\alpha$}
{\protect{\label{fig1}}}%
\end{figure}



\setlength{\unitlength}{1mm}%
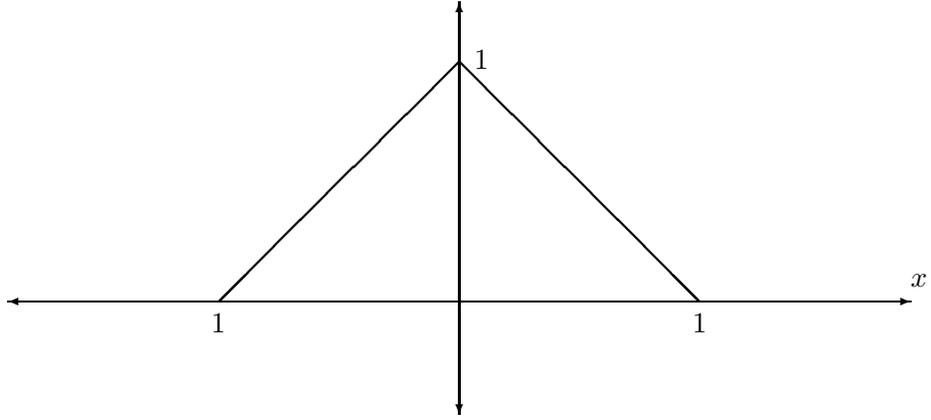
\begin{figure}
\begin{center}
\begin{picture}(140,45)(-68,-10)
\put(0,0){\vector(1,0){60}}
\put(0,0){\vector(-1,0){60}}
\put(0,0){\vector(0,1){40}}
\put(0,0){\vector(0,-1){15}}
\thicklines
\put(0,32){\line(1,-1){32}}
\put(0,32){\line(-1,-1){32}}


\put(60,2){$x$}
\put(31,-4){$1$}
\put(-33,-4){$1$}
\put(2,31){$1$}

\end{picture}
\end{center}
\caption{Graph of $\beta$}
{\protect{\label{fig2}}}%
\end{figure}



\bigskip

Instead of the unit grid we may consider a grid of squares
with sides of an arbitrary length $\rho >0$. Denote such a grid by ${\cal
G}_t^\rho$ if $t\in\R^2$ is a vertex of one of its members. Henceforth we
will call $\rho$ the {\it size} of the grid and $t$ its {\it reference point}.
We obtain another family of operators, defined by
$$
\PP_t^\rho f := \sum_{Q\in {\cal G}_t^\rho}
\sk{f}{h_{Q_0}}h_{Q_0}\,.
$$
Applying Proposition \ref{F} or modifying its proof, we can show the following.
\begin{prop}
\label{ro}
Choose $\rho>0$. Then averaging operators $\PP_t^\rho$ 
returns a convolution operator with the kernel
$$
F^\rho(x_1,x_2):=\frac{1}{\rho^2}\,F\left ( \frac{x_1}{\rho},
\frac{x_2}{\rho}\right )\,.
$$
\end{prop}

\noindent
Thus we have found the kernel of the operator, resulting
from averaging over {\sl all} grids of a {\sl fixed} size.
Our next step will be to average over all {\sl
sizes}. Let us explain what we mean by that.

Take $r>0$. A {\it lattice of calibre r} is
said to be a family of intervals (squares), obtained
from the standard dyadic lattice ${\cal L}$ by
dilating it by a factor $r$ and translating by an arbitrary
vector $t$. In other words, such a lattice (call it ${\cal L}_t^r$) is
the union of {\sl grids} of sizes $r\cdot 2^n\, ,n\in\Z$,
having $t$ as their reference point.


We introduce kernels
\begin{equation}
\label{k^r}
k^r:=\sum_{n=-\infty}^{\infty}F^{r\cdot 2^n}\,.
\end{equation}
By Proposition \ref{ro},
$$
\aligned
k^r(x)
& =\frac{1}{r^2}\sum_{n=-\infty}^{\infty}\frac{1}{4^n}\,F\Big(\frac{x}{r\cdot 2^n}\Big)\\
& =\frac{1}{r^2}\sum_{n=0}^{\infty}\frac{1}{4^n}\,F\Big(\frac{x}{r\cdot 2^n}\Big)+\frac{1}{r^2}\sum_{n=1}^{\infty}4^nF\Big(\frac{2^nx}{r}\Big)\,.
\endaligned
$$
The first sum in the lower line converges absolutely and uniformly on $\C$ because $F$ is bounded. If we restrict $x$ to the complement of any neighbourhood of the origin, then the second sum is actually ``uniformly" finite, because $F$ has a compact support. Consequently, the sum defining $k^r$ converges absolutely and uniformly on any such set.

The Proposition \ref{ro} provides also the relation
\begin{equation}
\label{diad.dilat}
k^r(2^mx)=\frac{k^r(x)}{4^m}\,
\end{equation}
for all $m\in\Z$, $r>0$, $x\in\C$.

The fact that $k^r*$ is a sum of operators, obtained by
averaging over grids of size $r\cdot 2^n$, hints at $k^r*$
itself being an average, this time over {\sl unions} of
these grids, i.e. lattices of calibre $r$. While it is not clear what could be a probability space corresponding to all lattices of a fixed calibre, we define the above-said average as a limit of averages of truncated lattices. Then the statement makes sense and holds \cite{DV}. Virtually the same proof establishes the lemma which follows below.

For $M\in\Z$ let the $M$-th partial sum of the series $k^r$ be
$$
k^r_M:=\sum_{n=-\infty}^{M}F^{r\cdot 2^n}\,.
$$
\begin{lema}
\label{Canon}
Function $k^r_M$ defines a bounded convolution operator on $L^p$. The
limit
$k^r*:=\lim_{M\rightarrow\infty} k_M^r*$
\noindent
exists in the strong sense and also gives rise to a bounded operator on
$L^p$.
\end{lema}

Next step is to average over dilations, in other words, over all calibres $r$. It is clear that the set of all possible calibres naturally corresponds to the
interval $[1,2)$. For our purpose, the most appropriate 
measure on this interval turns out to be $dr/r$. This
makes all other possible choices of intervals, e.g.
$[2^n,2^{n+1})$, have the same measure ($\log 2$). 

Averaging operators $k^r*$, i.e. integrating
$k^r$ with respect to the normalized measure $dr/r$, 
gives us a convolution operator once again. Call its kernel $k$. Then
$$
\aligned
k(x)& =\frac{1}{\log 2}\int_1^2k^r(x)\,\frac{dr}{r}
=\frac{1}{\log 2}\int_1^2\sum_{n=-\infty}^{\infty}F^{r\cdot
2^n}(x)\,\frac{dr}{r}\\
& =\frac{1}{\log 2}\sum_{n=-\infty}^{\infty}\int_1^2F^{r\cdot
2^n}(x)\,\frac{dr}{r}
=\frac{1}{\log 2}\sum_{n=-\infty}^{\infty}\int_{2^n}^{2^{n+1}}F^{s}
(x)\,\frac{ds}{s}
\endaligned
$$
\begin{equation}
\label{kk}
\hskip -154pt 
=\frac{1}{\log 2}\int_0^{\infty}F^{s}(x)\,\frac{ds}{s}
\,.
\end{equation}
Note that
\begin{equation}
\label{k}
k(x) = \frac{k
(x/|x|)} {|x|^2}
\end{equation}
for all $x\in\R^2\backslash\{0\}$. 
Because of this it suffices to know the behaviour of $k$ on $S^1$. By applying \eqref{kk} and Proposition \ref{ro} we get
\begin{equation}
\label{k|S}
k(e^{i\varphi})=\frac{1}{\log 2}\int_0^{\infty}F(re^{i\varphi})\,r\,dr\,.
\end{equation}

\bigskip
\noindent
{\bf Remark} 
This homogeneity is essential if we want to obtain the kernel of the Ahlfors-Beurling operator, which is, of course, homogeneous of the same degree. This, however, poses a restriction on the choice of the martingale transforms we are averaging. Namely, the coefficients corresponding to $h_{Q_*}$ should not depend on the size of $Q$. That is, they should be equal for all subgrids of a lattice. Yet in other words, the coefficients standing at $F^{r\cdot 2^n}$ in \eqref{k^r} should be the same for all $n$, otherwise we lose the homogeneity.

\bigskip

Finally, we are going to perform averaging over rotations. The rest of the proof will be somewhat different from the one, presented in \cite{DV}.

Choose $\psi\in [0,2\pi)$. Operators will be the same as before, just that the grids and lattices will consist of squares, rotated by the angle $\psi$ counterclockwise with respect to the standard position. Let $U_\psi: \C\rightarrow \C$ be defined by $U_\psi(\zeta):=\zeta e^{-i\psi}$. Then the convolution kernel of the operator $K_\psi$, which corresponds to the average over rotated lattices, is equal to $k_\psi:=k\circ U_\psi$. 
The operator itself satisfies the similarity relation
$
K_\psi=S_\psi^{-1}K_0S_\psi\,,
$
where $S_\psi f=f\circ U_{-\psi}$. 
\label{Tprime}

Both $k_\psi$ and $K_\psi$ appeared in \cite{DV}.

Now let us define a (weighted) average of operators $K_\psi$, which we denote by $T'$.
$$
\aligned
(T'f)(z):\!&=-\frac{1}{2\pi}\int_0^{2\pi}(K_\psi f)(z)\, e^{-2i\psi}\,d\psi\\
& =-\frac{1}{2\pi}\int_0^{2\pi}(k_\psi*f)(z)\, e^{-2i\psi}\,d\psi\\
& =-\frac{1}{2\pi}\int_0^{2\pi}
\int_\C k(\zeta e^{-i\psi})f(z-\zeta)\,dA(\zeta)
\, e^{-2i\psi}\,d\psi\,.
\endaligned
$$
Using the observation \eqref{k} we continue as
$$
\aligned
(T'f)(z)& =-\frac{1}{2\pi}\int_0^{2\pi}
\int_\C \frac{k(e^{i(\arg\zeta -\psi)})}{|\zeta|^2}\,f(z-\zeta)\,dA(\zeta)
\, e^{-2i\psi}\,d\psi\\
& =-\frac{1}{\pi}\int_\C\frac{f(z-\zeta)}{|\zeta|^2}\,
\frac{1}{2}\int_0^{2\pi}k(e^{i(\arg\zeta -\psi)})\, e^{-2i\psi}\,d\psi
\,dA(\zeta)\,.
\endaligned
$$
Now
$$
\int_0^{2\pi}k(e^{i(\arg\zeta -\psi)})\, e^{-2i\psi}\,d\psi=
e^{-2i\arg\zeta}\int_0^{2\pi}k(e^{i\varphi})\, e^{2i\varphi}\,d\varphi\,.
$$
From \eqref{k|S} it follows that \label{kFour}
$$
\aligned
\int_0^{2\pi}k(e^{i\varphi})\, e^{2i\varphi}\,d\varphi
& 
=\frac{1}{\log2}
\int_0^{2\pi}\int_0^{\infty}
F(re^{i\varphi})\,r\,dr
\, e^{2i\varphi}\,d\varphi\\
& =\frac{1}{\log2}
\int_\R\int_\R
F(x,y)\,\frac{(x+iy)^2}{x^2+y^2}
\,dx\,dy
\,,
\endaligned
$$
hence we arrived at the final result.
\begin{thm} 
\label{C}
Assuming the notation as above,
$$
T=C\cdot T'\,,
$$
where
\begin{equation}
\label{konst}
\frac{1}{C}
=\frac{1}{2\log2}
\int_\R\int_\R
F(x,y)\,\frac{(x+iy)^2}{x^2+y^2}
\,dx\,dy
\,.
\end{equation}
\end{thm}

\bigskip
\noindent
Although this theorem is, at the first glance, just a worked-out version of \eqref{grobo}, it deserves to be written separately. For it gives us a precise formula for the constant arising from the 
averaging process. Moreover, it can be applied to many different Haar bases.
Once fixing the Haar system, the kernel $F$ depends only on the coefficients of the ``core" operator. 
For example, if 
$$
\PP_t f := \sum_{Q\in {\cal G}_t} [\sigma^1\sk{f}{h_{Q}^1}h_{Q}^1+\sigma^2\sk{f}{h_{Q}^2}h_{Q}^2+\sigma^3\sk{f}{h_{Q}^3}h_{Q}^3]\,,
$$
for certain $\sigma^1$, $\sigma^2$, $\sigma^3\in S^1$, then, according to Proposition \ref{F}, 
$$
F(x,y)=-\sigma^1\beta(x_1)\alpha(x_2)-\sigma^2\alpha(x_1)\beta(x_2)+\sigma^3\alpha(x_1)\alpha(x_2)\,.
$$
Thus looking for the best constant, corresponding to the chosen Haar system, comes down to maximizing the integral in \eqref{konst} for all admissible combinations of $\sigma$'s.


The constant $C$, arising from ${\cal H}_{orig}$, is too big. Our way of attempting its decrease is to consider more general Haar systems. However, we should restrict ourselves to searching among systems for which we are able to obtain $p^*-1$ estimates (and ${\cal H}_{orig}$ is not one of them).

\subsection{Optimizing coefficients}

We saw that the idea to change the ``original" Haar system emerged from our inability to prove $p^*-1$ estimates for it. Now that we introduced appropriate changes in Section \ref{remedy}, we are going to repeat the averaging process and use Theorem \ref{C} to calculate the best constant for the modified system. We will gradually extend our search to more general systems including those assuming different geometrical properties.

Let us start by considering ${\cal H}_{new}$, defined in Section \ref{remedy}. 
Choose complex numbers $\sigma_0$, $\sigma_+$ and $\sigma_-$ with modulus one. Let us examine the operators of the type
\begin{equation}
\label{P_t}
{\cal P}_t f:=
\sum_{Q\in {\cal G}_t}
\left[\sigma_{0}\sk{f}{h_{Q_0}}h_{Q_0}
+\sigma_{+}\sk{f}{h_{Q_+}}h_{Q_+}
+\sigma_{-}\sk{f}{h_{Q_-}}h_{Q_-}\right].
\end{equation}
The coefficients 
$\sigma_0$, $\sigma_+$, $\sigma_-$ 
are chosen not to depend on squares $Q$, for otherwise we might already get in trouble when trying to run the first averaging process -- the one over translations. 

It is convenient to write the 
summands in terms of the functions from ${\cal H}_{orig}$, since for them the kernels resulting after the averaging were already computed. By using the identities \eqref{new} from page \pageref{new} we get 
$$
\sigma_{0}H_Q^1f
+\frac{\sigma_{+}+\sigma_{-}}{2}\,(H_Q^{2}+H_Q^{3})f
+\frac{\sigma_{+}-\sigma_{-}}{2}\,(\sk{f}{h_Q^{2}}h_Q^{3}+\sk{f}{h_Q^{3}}h_Q^{2})
$$
where $H_Q^{j}f=\sk{f}{h_Q^{j}}h_Q^{j}$, $j=1,2,3$.

Let us average operators ${\cal P}_t$ over grids. 
A proof, analogous to that of Proposition \ref{F}, shows that the sum of mixed terms in parentheses on the right becomes zero. Thus, by Proposition \ref{F}, the kernel we get is 
$$
F(x,y)=-\sigma_{0}\beta(x)\alpha(y)+\frac{\sigma_{+}+\sigma_{-}}{2}\,(-\alpha(x)\beta(y)+\alpha(x)\alpha(y))\,.
$$
We can assume that $\sigma_{0}=1$, for we are only interested in the maximum of the absolute value of \eqref{konst}. Since the integral of $-\alpha(x)\beta(y)+\alpha(x)\alpha(y)$ will, of course, be real, the maximum will be obtained when $\sigma_{+}=\sigma_{-}=1$ or $\sigma_{+}=\sigma_{-}=-1$. The first choice would mean that we are eventually averaging the identity operators, thus it has to be discarded. \label{sigma} This is how we obtained the best coefficients $\sigma$ in the case of ${\cal H}_{new}$. Hence our problem boils down to estimating one simple concrete martingale transform. We saw that our best bet is
\begin{equation}
\label{DSCH}
F(x,y)=\alpha(x)\beta(y)-\beta(x)\alpha(y)-\alpha(x)\alpha(y)
\end{equation}
and the evaluation of the constant $C$ then returns approximately 2,07 (the exact value is given on page \pageref{izraz}).

\subsection{Parallelograms}

We can perform the same averaging process for more general Haar systems. One way of introducing parameters of generality is to consider functions supported on parallelograms rather than squares. 

Let us start with a parallelogram whose 
sides have lengths 1 and $b$, respectively, 
and whose inner angle at the bottom left corner equals $\varphi$ (see Figure \ref{paral}). Here $b>0$ and $\varphi\in (0,\pi)$ can be arbitrary. 
In addition, we can assume that its lower side is the segment $[0,1)\times \{0\}$.
We cover $\R^2$ by a grid of similar parallelograms and consecutively form a corresponding dyadic lattice. To each of its members we assign, as always so far, three Haar functions. Figure \ref{paral} displays the analogue of $h_{Q_0}$ (for which we retain the same name) for the ``generic" parallelogram. It is obvious how the other two functions, $h_{Q_+}$ and $h_{Q_-}$, should look like.
The set of all Haar functions $h_{Q_0}$, $h_{Q_+}$ and $h_{Q_-}$, where $Q$ runs over the above-said dyadic lattice, will be denoted by ${\cal H}_{b,\varphi}$. In particular, ${\cal H}_{new}={\cal H}_{1,\pi/2}$. This construction was clearly made to fit the proof of Theorem \ref{}, in other words, 
$$
\begin{array}{l}
\text{\it the corresponding martingale transforms also admit } L^p-\text{\it norms}\\
\text{\it smaller than } p^*-1\,.
\end{array}
$$


\setlength{\unitlength}{1mm}
\begin{figure}
\begin{center}
\begin{picture}(140,-10)(-55,0)


\put(0,-2){\line(0,-1){3}}
\put(0,-3.5){\vector(1,0){16}}
\put(16,-2){\line(0,-1){3}}
\put(16,-3.5){\vector(-1,0){16}}

\put(7,-7){$1$}

\put(18,-1){\line(2,-1){4}}
\put(20,-2){\vector(1,2){12.5}}
\put(30.5,24){\line(2,-1){4}}
\put(32.5,23){\vector(-1,-2){12.5}}

\put(27,7){$b$}
\put(16.25,17.5){$+$}
\put(10,5){$-$}

 \put(1.9,0){\oval(8,8)[tr]}
 \put(1.6,1.3){$\varphi$}


\multiput(6.25,12.5)(3.5,0){5}{\line(1,0){2.2}}
\put(6.25,12.5){\line(1,0){16}} 
\thicklines
\put(0,0){\line(1,2){12.5}}
\put(16,0){\line(1,2){12.5}}
\put(0,0){\line(1,0){16}}
\put(12.5,25){\line(1,0){16}}


\end{picture}
\end{center}
\caption{General $h_{Q_0}$}
{\protect{\label{paral}}}%
\end{figure}
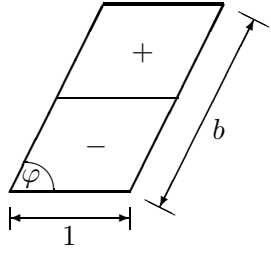


So let us define ${\cal P}_t$ as in \eqref{P_t}, just that this time the grid ${\cal G}_t$ consists of parallelograms. A similar consideration as before shows that the choice $\sigma_0=1$, $\sigma_+=\sigma_-=-1$ is optimal. In this setting, too, we can use 
Theorem \ref{C}, i.e. we can represent $T$ as an average of martingale transforms, arising from ${\cal P}_t$, over 
all lattices consisting of parallelograms of sides $2^nb$, ¢$n\in\Z$, and inclination $\varphi$. The r\^ole of $F$ is now assumed by a more general kernel $F_{b,\varphi}$, given by the formula
\begin{equation}
\label{F_b}
F_{b,\varphi}(x,y)=\frac1{b\sin\varphi}\, F\Big(x-\frac{y}{\tan\varphi},\frac{y}{b\sin\varphi}\Big)\,.
\end{equation}
Also the formula for $C$ now depends on two parameters, explicitly, it takes the form 
$$
C=\frac{2\log 2}{I(b,\varphi)}\,,
$$
where
$$
\aligned
I(b,\varphi) & = \int_\R\int_\R F_{b,\varphi}(x,y)\frac{x+iy}{x-iy}\,dx\,dy \\
& =-2(b\sin\varphi)^2\int_\R\int_\R F(x,y)\,\frac{y^2}{x^2+(2b\cos\varphi) xy +b^2y^2}\,dx\,dy\,. \\
\endaligned
$$
Recall that $F$ is given by \eqref{DSCH}.

Numerical tests indicate that the maximum of all $I(b,\varphi)$ is probably attained when $\varphi=\pi/2$ and $b=\sqrt{2}$, in other words, for rectangles with ratio of sides equal to $\sqrt{2}$. Note that such rectangles have the peculiar property --unique among all parallelograms-- that cutting them in half (the way we do while constructing our Haar system) 
preserves the ratio between the longer and the shorter side as well as the ``inclination".
Though, we are still not certain on how this geometric property relates 
to the fact that such rectangles yield the best constant among all parallelograms.

Let us compute its exact value.
$$
\aligned
I(\sqrt{2},\pi/2) & = 
-4\int_{-1}^1\int_{-1}^1 [\alpha(x)\beta(y)-\beta(x)\alpha(y)-\alpha(x)\alpha(y)] 
\frac{y^2}{x^2+2y^2}\,dx\,dy\\
& =-16\int_0^1\int_0^1 [\alpha(x)\beta(y)-\beta(x)\alpha(y)-\alpha(x)\alpha(y)] 
\frac{y^2}{x^2+2y^2}\,dx\,dy\\
& = 16\int_0^{1/2}\int_0^{1/2}(9xy-5x-y+1)\,\frac{y^2}{x^2+2y^2}\,dx
\ dy\\
& \hskip 15pt +16\int_0^{1/2}\int_{1/2}^1 (1-x)(7y-3)\,\frac{y^2}{x^2+2y^2}\,dx\ dy\\
& \hskip 15pt +16\int_{1/2}^1
\int_{0}^1 (1-x)(1-y)\,\frac{y^2}{x^2+2y^2}\,dx\ dy\,.
\endaligned
$$
Since
$$
\int_0^{1/2}\frac{9xy-5x-y+1}{x^2+2y^2}\,dx=\frac{9y-5}{2}\log\left(1+\frac{1}{8y^2}\right)+\frac{1-y}{y\sqrt{2}}\arctan\frac{1}{2y\sqrt{2}}\,,
$$
the first integral equals\\

$
\displaystyle{16\int_0^{1/2}
y^2
\left(
\frac{9y-5}{2}\log\left(1+\frac{1}{8y^2}\right)+\frac{1-y}{y\sqrt{2}}\arctan\frac{1}{2y\sqrt{2}}
\right)
\, dy}
$
$$
\aligned
=\frac{1}{96}
\Big(
& -26+40\sqrt{2}\,\pi -16\sqrt{2}\arctan(1/\sqrt{2})-48\sqrt{2}\arctan\sqrt{2}\\
& +52\log 2-71\log 3
\Big)\,.
\endaligned
$$

Similarly, the calculation of the second integral gives
$$ 
\frac{1}{96}\left(
-14+256\sqrt{2}\arctan(1/\sqrt{2})-160\sqrt{2}\arctan\sqrt{2}+124\log 2-77\log 3
\right)\,,
$$
while the third one turns out to be 
$$ 
\frac{1}{24}\left(
10+48\sqrt{2}\arctan(1/\sqrt{2})-32\sqrt{2}\arctan\sqrt{2}+20\log 2-11\log 3
\right)\,.
$$
Thus altogether we get
$$ 
I(\sqrt{2},\pi/2) =
\frac{1}{\sqrt{2}}\left(
\frac{5\pi}{6}+9\arctan\frac{1}{\sqrt{2}}-7\arctan\sqrt{2}
\right)
+\frac{8}{3}\log 2-2\log 3\,.
$$
\begin{cor}
\label{007}
For every $p\in (1,\infty)$,
\begin{equation}
\label{2007}
\frac{\nor{T}_{p}}{p^*-1}\leqslant\frac{2\log 2}{I(\sqrt{2},\pi/2)}\approx 2,00714\,.
\end{equation}
\end{cor}

The preceding analysis showed that the optimal results were obtained for coefficients $\sigma$ which were not only independent on the position of the Haar function (i.e. its support $Q$), but also assumed real values, i.e. $\pm 1$ (see p. \pageref{sigma}). \label{concrete}
Thus if for a fixed lattice $\mathcal{L}$ we denote by $P_*$ the projection onto the closed linear subspace generated by all Haar functions $\mn{h_{Q_*}}{Q\in\mathcal{L}}$, then obviously 
$$
\ \ P_0+P_++P_-
$$
is the identity operator, while by slightly changing it to
$$
-P_0+P_++P_-
$$
we obtain exactly those martingale transforms that yielded the optimal constant. These may be viewed as the simplest nontrivial martingale transforms.  (We should emphasize that ``nontriviality" here necessarily assumes a different meaning than in the one-dimensional case of Burkholder's theorem. There only one Haar function is assigned to only interval in the lattice, which forces nontrivial martingale transforms' coefficients to depend on the intervals. In higher dimensions, as noted above, this is not the case.) 
Hence it is to be expected that their norm estimate is much better than $p^*-1$, 
which is the supremum of norms over all possible martingale transforms $T_\sigma$. 
Though we do not know how to derive bounds for them. Possibly that would require a skillful improvement of the Bellman function of Burkholder by using the additional information supplied by these concrete and simple operators. 

\bigskip
\noindent
{\bf Remark}
There is another reason which hinders us from proving \eqref{domneva} 
in its entirety by means of the method described in this work. Namely, whenever we estimate a norm of the operator resulting from some averaging process, we lose a bit due to the triangle inequality. For example, the outcome of averaging over translations is simply estimated as   $\nor{\int\PP_t\,dt}\leqslant\int\nor{\PP_t}\,dt$ and similarly for dilations and rotations. We do not have any control over how much we actually lose.


\subsection{Higher dimensions}
\label{hidi}

A cutting, analogous to that of Section \ref{remedy}, and consequent estimates can be successfully tried in arbitrary $\R^n$. For example, the case of $\R^3$ is pretty straightforward:

Take a cube $Q$ and denote by $Q_1$, $Q_2$, $\dots, Q_8$ some ordering of its dyadic subcubes. 
Let $\chi_{j,k}$ stand for $\chi_{Q_j}+\chi_{Q_k}$. To $Q$ we can assign 7 Haar functions, namely
$$
\aligned 
& \chi_{1,2,3,4}-\chi_{5,6,7,8}  \\
& \chi_{1,2}-\chi_{3,4}  \\
& \chi_{5,6}-\chi_{7,8}  \\
& \chi_1-\chi_2 \\
& \chi_3-\chi_4 \\
& \chi_5-\chi_6\\
& \chi_7-\chi_8 
\endaligned
$$
plus the normalization in $L^2(\R^3)$. If we let $Q$ run over some dyadic lattice in $\R^3$, we get a complete orthonormal system in $L^2(\R^3)$. Applying the same proof as in $\R^2$ shows that the associated martingale transforms are bounded by $p^*-1$ on $L^p(\R^3)$. 

In a completely analogous way we proceed in arbitrary dimension. The Haar system in $\R^n$ constructed in this manner will be denoted by ${\cal H}_{new}^n$. The nature of construction gives rise to the following result, whose proof we omit due to its similarity to the proof of Theorem \ref{}.

%

\begin{thm}
\label{mart-n}
Let $T_\sigma$ be any martingale transform arising from ${\cal H}_{new}^n$. Then $\nor{T_\sigma}_p\leqslant p^*-1$ and this estimate is sharp.
\end{thm}

We can repeat the averaging of such martingale transforms, as done on page \pageref{Tprime}. As a result we get a kernel $k$ --an exact analogue of the one in \eqref{k}-- with the property
$$
k(x)=\frac{k(x/|x|)}{|x|^n}\,
$$
for $x\in\R^n\backslash\{0\}$. Certainly, the norm of the operator $f\mapsto k*f$ is majorized by $p^*-1$ on $L^p(\R^n)$. 
Now take a convolution operator $S$ on $L^p(\R^n)$ whose kernel $\tau$ is continuous on $S^{n-1}$ and admits the same homogeneity as $k$, namely
\begin{equation}
\label{homog}
\tau(x)=\frac{\tau(x/|x|)}{|x|^n}\,.
\end{equation}
We would like to find some sufficient conditions on $\tau$ which would allow a representation of $S$ as an average of martingale transforms.

So far we treated the case $n=2$ and $\tau(e^{i\psi})=e^{-2i\psi}$. Its advantage was that we could identify the sphere $S^1$ with the group of rotations. 

Let $G$ be $SO(n)$, the group of rotations in $\R^n$. Its Haar measure will be denoted by $dU$. For $U\in G$ write $k_U=k\circ U$ and $K_Uf=k_U*f$, as before. Then $K_U$ is the operator occurring after averaging over all dyadic lattices, rotated by $U$. Instead of  $e^{-2i\psi}$ now we have the weight $w(U)\in L^1(G)$. We repeat the proof from page \pageref{C}. Define, for $x\in\R^n$,
$$
\aligned
(T'f)(x):\!&=
\int_G(K_U f)(x)\, w(U)\,dU\\
& =
\int_G((k\circ U)*f)(x)\, w(U)\,dU\\
& =
\int_G\int_{\R^n} k(Uy)f(x-y)\,dy
\, w(U)\,dU\\
& =
\int_G\int_{\R^n} \frac{k(Uy/|y|)}{|y|^n}\,f(x-y)\,dy
\, w(U)\,dU\\
& =\int_{\R^n}\frac{f(x-y)}{|y|^n}\,
\int_Gk(Uy/|y|)\, w(U)\,dU\,dy\,.
\endaligned
$$
Thus in order to conclude that $T'=S$,
we would like to have 
\begin{equation}
\label{tauw}
\int_Gk(U\xi)\, w(U)\,dU=
\tau(\xi)\,,\hskip 20pt \xi\in S^{n-1}\,.
\end{equation}
Since from the construction of $T'$ it follows that $\nor{T'}_p\leqslant \nor{w}_{L^1(G)}(p^*-1)$, that would obviously imply the same bound for $S$. For that reason we are going to take \eqref{tauw} as a {\sl definition} of $\tau$ (depending on $w$). Let us make the statement more precise. 
\begin{thm}
Take a function $w\in L^1(G)$ and define $\tau\in C(\R^n\backslash\{0\})$ by  
$$
\tau(\xi)=\int_Gk(U\xi)\, w(U)\,dU\,
$$
for $\xi\in S^{n-1}$ and by \eqref{homog} on $\R^n\backslash\{0\}$. Then the operator $S$ of convolution with $\tau$ admits on $L^p(\R^n)$ the estimate
$$
\nor{S}_p\leqslant \nor{w}_{L^1(G)}(p^*-1)
$$
for every $p\in (1,\infty)$.
\end{thm}


\section{Appendix: averaging in other systems}
\label{largo}

\subsection{Horizontal vs. diagonal}
\label{hvd}

The Haar functions in ${\cal H}_{new}$ and their subsequent generalizations in ${\cal H}_{b,\varphi}$ were defined in the following way. On each $Q$ we took the first Haar function to be $h_Q^1$. Then we cut $Q$ in two horizontally, in other words, $h_{Q_+}$ and $h_{Q_-}$ were supported, respectively, on preimages of $(0,\infty)$ and $(-\infty,0)$ under $h_Q^1$. Of course, such cutting can be performed starting with functions other than $h_Q^1$. If we begin with $h_Q^2$, we do not encounter any novelties. This is, roughly, because $h_Q^2$ can be viewed as a rotation of $h_Q^1$ by $\pi/2$. Since we also average over rotations, it comes to no surprise that this case yields the same best constant. 

But the case of $h_Q^3$ is different, in that its preimages of the positive and negative half-line are geometrically distinct from what we had before. Let us take a closer look of what this observation implies.

Firstly, let us introduce some notation. In order to avoid confusion, we will denote the Haar functions for such a system by $h_{Q_0}^d$, $h_{Q_+}^d$ and $h_{Q_-}^d$ (here $d$ stands for ``diagonal"). In case of squares --and analogously for parallelograms--, they are given by  
%
%
%
%
%
$$ 
h_{Q_0}^d \equiv
\begin{array}{|c|c|}
\hline \raisebox{0pt}[12pt][6pt]{$-$} & +  \\ \hline
\raisebox{0pt}[12pt][6pt]{$+$} &  -  \\ \hline
\end{array}
\hskip 30pt 
h_{Q_+}^d \equiv
\begin{array}{|c|c|}
\hline \raisebox{0pt}[12pt][6pt]{} & +  \\ \hline
\raisebox{0pt}[12pt][6pt]{$-$} &    \\ \hline
\end{array}
\hskip 30pt 
h_{Q_-}^d \equiv
\begin{array}{|c|c|}
\hline \raisebox{0pt}[12pt][6pt]{$+$} &   \\ \hline
\raisebox{0pt}[12pt][6pt]{} &  -  \\ \hline
\end{array}
$$
The set of all such Haar functions, where $Q$'s run over the lattice of parallelograms of type $(b,\varphi)$, will be denoted by ${\cal H}_{b,\varphi}^d$. Once again we define ${\cal P}_t$ as in \eqref{P_t}. As before, we can assume that $\sigma_0=1$. Then
$$
{\cal P}_tf=\sum_{Q\in{\cal G}_t}[H_Q^3f
+\frac{\sigma_{+}+\sigma_{-}}{2}(H_Q^{1}+H_Q^{2})f
+\frac{\sigma_{+}-\sigma_{-}}{2}(\sk{f}{h_Q^{1}}h_Q^{2}+\sk{f}{h_Q^{2}}h_Q^{1})]\,.
$$ 
This time, however, averaging the mixed term $\sk{f}{h_Q^{1}}h_Q^{2}+\sk{f}{h_Q^{2}}h_Q^{1}$ over translations does {\sl not} return zero. Instead, we get a convolution operator whose kernel is $-2\gamma(x)\gamma(y)$, in case of squares, where $\gamma=h_0*\chi_0$ (notation as in the formulation of Proposition \ref{F}). The graph of $\gamma$ is given in Figure \ref{gama}.


\setlength{\unitlength}{1mm}
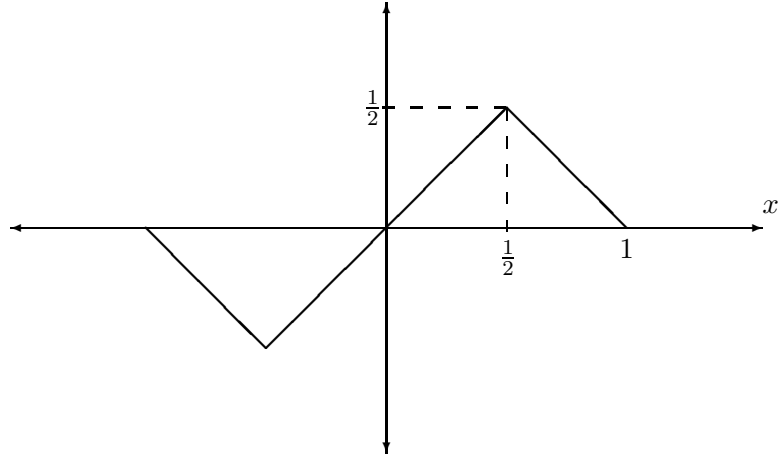
\begin{figure}
\begin{center}
\begin{picture}(140,60)(-68,-25)
\put(0,0){\vector(1,0){50}}
\put(0,0){\vector(-1,0){50}}
\put(0,0){\vector(0,1){30}}
\put(0,0){\vector(0,-1){30}}
\multiput(-0.5,16)(3.72,0){5}{\line(1,0){1.5}}
\multiput(16,-0.5)(0,3.7){5}{\line(0,1){1.5}}

\put(15,-5){$\frac12$}
\put(-3,15){$\frac12$}

\put(50,2){$x$}
\put(31,-4){$1$}
\put(-12,16){\line(1,-1){2}}

\thicklines
\put(-16,-16){\line(1,1){32}}
\put(16,16){\line(1,-1){16}}
\put(-32,0){\line(1,-1){16}}


\end{picture}
\end{center}
\caption{Graph of $\gamma$}
{\protect{\label{gama}}}%
\end{figure}


Hence the kernel of averages of ${\cal P}_t$ over rotations is 
$$
F(x,y)=\alpha(x)\alpha(y)-\frac{\sigma_{+}+\sigma_{-}}{2}[\alpha(x)\beta(y)+\beta(x)\alpha(y)]-(\sigma_+-\sigma_-)\gamma(x)\gamma(y)\,.
$$
Note that earlier $F$ was even both in $x$ and $y$. Now this is not the case anymore. 
Furthermore, when we are averaging over grids of {\bf parallelograms} with sides 1 and $b$ and inclination $\varphi$ (see Figure \ref{paral}), the kernel is  
$$
\frac1{b\sin\varphi} F\Big(x-\frac{y}{\tan\varphi},\frac{y}{b\sin\varphi}\Big)\,.
$$
Let us 
restrict ourselves to {\bf rectangles}, i.e. choose $\varphi=\pi/2$. Then the constant can be evaluated from the formula \eqref{konst}:
$$
\aligned
\frac{1}{C} & =\frac{1}{2\log2}
\int_\R\int_\R
\frac1b F\Big(x,\frac{y}{b}\Big)\,\frac{(x+iy)^2}{x^2+y^2}
\,dx\,dy\\
& =\frac{1}{2\log2}
\int_\R\int_\R
F(x,t)\,\frac{(x+i\,bt)^2}{x^2+b^2t^2}
\,dx\,dt
\,.
\endaligned
$$
So far the kernels $F$ were even with respect to both variables, so we only had to consider the real part of the weight in the integral above. This time it is not so, which can be seen as a good news, since the kernel of $T$ itself is complex, i.e. non-real. We split $F$ into a sum of two functions, one being even and the other odd. That is, $F=F_1-F_2$, where
$$
\aligned
F_1(x,y)&=\alpha(x)\alpha(y)-\frac{\sigma_{+}+\sigma_{-}}{2}[\alpha(x)\beta(y)+\beta(x)\alpha(y)]\\
F_2(x,y)&=(\sigma_+-\sigma_-)\gamma(x)\gamma(y)\,.
\endaligned
$$
Consequently,
\begin{equation}
\label{kajvem}
\aligned
\frac1C=& \frac{1}{2\log2}
\int_\R\int_\R
F_1(x,t)\,\frac{x^2-b^2t^2}{x^2+b^2t^2}
\,dx\,dt\\
& -i\,\frac{b}{\log2}
\int_\R\int_\R
F_2(x,t)\,\frac{xt}{x^2+b^2t^2}
\,dx\,dt\,.
\endaligned
\end{equation}
By using the following simple observations 
\begin{itemize}
\item
$\alpha$ and $\beta$ are even, $\gamma$ is odd
\item
all three functions are supported on the interval $[-1,1]$ 
\item
$\int_0^1\alpha(x)\,dx=0$
\end{itemize}
and the equality
$$
\frac{x^2-b^2t^2}{x^2+b^2t^2}=2\,\frac{x^2}{x^2+b^2t^2}-1\,,
$$
the formula \eqref{kajvem} simplifies into
$$
\aligned
\frac1C=& \frac{4}{\log2}
\int_0^1\int_0^1
[\alpha(x)\alpha(y)-\frac{\sigma_{+}+\sigma_{-}}{2}(\alpha(x)\beta(y)+\beta(x)\alpha(y))]\,\frac{x^2}{x^2+b^2y^2}
\,dx\,dy\\
& -i\,\frac{4b}{\log2}(\sigma_+-\sigma_-)
\int_0^1\int_0^1
\gamma(x)\gamma(y)\,\frac{xy}{x^2+b^2y^2}
\,dx\,dy\,.
\endaligned
$$
It is easy to see that, for arbitrary $b$, the maximum of the expression above over all $\sigma_+,\sigma_-\in S^1$ occurs when $\sigma_-=\overline\sigma_+$. Write $\sigma_\pm=\cos\vartheta\pm i\sin\vartheta$. Then  
$$
\aligned
\frac{\log 2}{4C}=& 
\int_0^1\int_0^1
\alpha(x)\alpha(y)\,\frac{x^2}{x^2+b^2y^2}
\,dx\,dy\\
&-\cos\vartheta\int_0^1\int_0^1[\alpha(x)\beta(y)+\beta(x)\alpha(y)]\,\frac{x^2}{x^2+b^2y^2}
\,dx\,dy\\
& +2b\sin\vartheta
\int_0^1\int_0^1
\gamma(x)\gamma(y)\,\frac{xy}{x^2+b^2y^2}
\,dx\,dy\,.
\endaligned
$$
Observe that 
\begin{itemize}
\item
$b\mapsto 1/b$ changes the sign of the 1. and 2. integral above, while the 3. remains unchanged
\item
$\vartheta\mapsto-\vartheta$ changes the sign of the 3. integral above, while the 1. and the 2. remain unchanged.
\end{itemize}
Thus if we write $C=C(b,\vartheta)$, we see that 
$C(1/b,\vartheta)=-C(b,-\vartheta)$. Hence we can restrict our search to maximizing over $b\in(0,1]$, $\vartheta\in(-\pi,\pi)$. 
However, preliminary numerical testing hints that the best constant for such Haar systems might still be worse than the one obtained in Corollary \ref{007}, after all. 


\subsection{Triangles}

We shall describe yet another Haar-type system. 
This time a grid, instead of consisting of squares, consists of triangles (see Figure \ref{figgrid}). A lattice is made --like in the ``rectangular" case-- of unions of grids, where the objects in the next grid are twice as large (i.e. four times in terms of area).


\setlength{\unitlength}{1mm}
\begin{figure}
\begin{center}
\begin{picture}(140,-10)(-70,0)
\thicklines
\put(-32,0){\line(1,0){64}}
\put(-32,16){\line(1,0){64}}
\put(-32,32){\line(1,0){64}}
\multiput(-32,0)(16,0){5}{\line(0,1){32}}
\multiput(-16,0)(16,0){4}{\line(-1,1){16}}
\multiput(-16,16)(16,0){4}{\line(-1,1){16}}
\put(-34,-2){$t$}


\end{picture}
\end{center}
\caption{Grid ${\cal G}_t$}
{\protect{\label{figgrid}}}%
\end{figure}
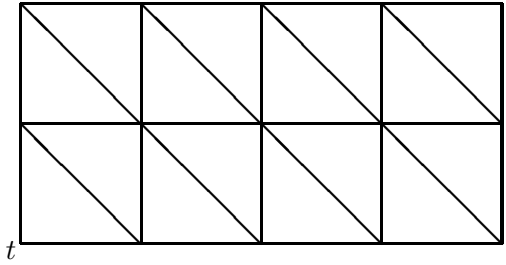


To each triangle from a grid we assign three Haar functions, as plotted on Figures \ref{fig0}, \ref{fig+} and \ref{fig-}. However, in a grid we have two kinds of triangles: those in ``normal" position and those ``turned upside down". But then again, Haar functions corresponding to the latter ones are easy to guess. 
Thus we have a complete orthonormal system with 3 different types of functions, namely $h_0$, $h_+$ and $h_-$.
The same procedure of obtaining sons (i.e. joining the middle points of triangle's sides) of equal shape certainly works for arbitrary triangles and not only for isosceles ones. 
The Haar system, whose building block is the triangle from Figure \ref{ab}, will be labeled as ${\cal H}_{a,b}^\vartriangle$. 
Hence we have a whole family of derived martingale transforms (denoted 
by $T_\sigma$, as usual). Of course, these Haar functions were constructed 
so that $T_\sigma$ should admit $p^*-1$ estimates on $L^p$.


\setlength{\unitlength}{1mm}%
\begin{figure}
\begin{center}
\begin{picture}(140,45)(-70,0)
\put(-5,0){\line(1,0){50}}
\put(0,-5){\line(0,1){45}}
\put(32,0){\line(1,2){10.8}}
\thicklines
\put(0,0){\vector(1,0){32}}
\put(0,0){\vector(2,1){42.8}}
\put(43,21){$(a,b)$}
\put(31,-4){$1$}
\put(-3,31){$1$}
\put(-1,31){--}


\end{picture}
\end{center}
\caption{Building block of ${\cal H}_{a,b}^\vartriangle$}
{\protect{\label{ab}}}%
\end{figure}


\subsubsection*{Calculating the kernel}

We proceed to the averaging. We will perform it starting with ${\cal H}_{a,b}^\vartriangle$, where $a,b\in\R^2$, $b\not =0$, are fixed. The notation ${\cal G}_t$ will stand for the grid, generated by the triangle from Figure \ref{ab}, translated by $t\in\R^2$. 
Introduce the operator 
$$
\PP_t: f\mapsto \sum_{D\in {\cal G}_t}\sk{f}{h_{D}}h_{D}\,.
$$
Our objective is computing the average 
$$
\frac{1}{|\Omega|}\int_{\Omega} \PP_t\,dt\,.
$$
The probability space (and domain of integration) $\Omega$ can be defined in the following way. Let $\Omega^+$ be a triangle from ${\cal G}_0$ which contains 0. 
Set $\Omega^-:=-\Omega^+$, which is again a triangle from ${\cal G}_0$ (modulo the boundary) but of opposite facing. Then we define $\Omega:=\Omega^+\cup\Omega^-$. 
%
%
\begin{prop}
\label{Piazzolla}
$$
\frac{1}{|\Omega|}\int_{\Omega} \PP_tf\,dt=F*f\,,
$$
where
$$
F(z)=\frac{2}{|\Omega|}\int_{\Omega^+} h_{\Omega^+}(\zeta)h_{\Omega^+}(\zeta-z)\,dA(\zeta)\,.
$$
\end{prop}

In order to avoid congestion of symbols and increase transparency in the proof we omitted adding $a,b$ to the symbol for the grid and also specifying the type of the Haar functions considered. The Proposition is valid for any type (of course, the one in the definition of $\PP_t$ should agree with the one appearing in the formula for $F$).\\ 

\dok
Our path will not differ much from the one followed in Proposition \ref{F}. 
For each triangle $D$ denote $H_Df=\sk{f}{h_{D}}h_{D}$. Fix $x\in\R^2$. Then 
$$
\PP_tf(x)=\sum_{D\in {\cal G}_t}H_Df(x)=\int_{\R^2}f(s)h_{D}(s)\,ds\,h_{D}(x)
$$
for the unique $D\in{\cal G}_t$ which contains $x$. Without any loss of generality we can assume that $D=\Omega_t^+:=\Omega^++t$ or $D=\Omega_t^-:=\Omega^-+t$. Then $x\in\Omega_t^+\cup\Omega_t^-$, which is equialent to $t\in(x-\Omega^+)\cup(x-\Omega^-)=x-\Omega\equiv\Omega$. The latter identification is to be understood in the sense that $x-\Omega$ and $\Omega$ are equal as probability spaces for our averaging process. Hence\\

$
\displaystyle{
\frac{1}{|\Omega|}\int_{\Omega}\sum_{D\in {\cal G}_t}H_Df(x)\,dt 
}
$
\vskip -5pt
$$
\aligned
 =\frac{1}{|\Omega|}\bigg( & 
\int_{x-\Omega^+}\int_{\R^2}f(s)h_{\Omega_t^+}(s)\,ds\,h_{\Omega_t^+}(x)\,dt\\
+& \left. 
\int_{x-\Omega^-}\int_{\R^2}f(s)h_{\Omega_t^-}(s)\,ds\,h_{\Omega_t^-}(x)\,dt
\right)\\
=\frac{1}{|\Omega|}\bigg( & 
\int_{x-\Omega^+}\int_{\R^2}f(s)h_{\Omega^+}(s-t)\,ds\,h_{\Omega^+}(x-t)\,dt\\
+ & \left. 
\int_{x-\Omega^-}\int_{\R^2}f(s)h_{\Omega^-}(s-t)\,ds\,h_{\Omega^-}(x-t)\,dt
\right)
\endaligned
$$
$$
\aligned
\hskip 4pt =\int_{\R^2}f(s)\frac{1}{|\Omega|}\bigg( &  
\int_{x-\Omega^+}h_{\Omega^+}(s-t)\,h_{\Omega^+}(x-t)\,dt\\
+ & 
\int_{x-\Omega^-}h_{\Omega^-}(s-t)\,h_{\Omega^-}(x-t)\,dt\bigg)\,ds\,.
\endaligned
$$
By entering $w:=x-t$, the first inner integral in the parentheses above can be rewritten as 
$$
\int_{\Omega^+}h_{\Omega^+}(s-x+w)\,h_{\Omega^+}(w)\,dw\,.
$$
Equally, the second integral is equal to 
$$
\aligned
\int_{\Omega^-}h_{\Omega^-}(s-x+w)\,h_{\Omega^-}(w)\,dw & \\
\stackrel{(\star 1)}{=} & \int_{\Omega^-}h_{\Omega^+}(x-s-w)\,h_{\Omega^+}(-w)\,dw\\
\stackrel{(\star 2)}{=} & \int_{\Omega^+}h_{\Omega^+}(x-s+v)\,h_{\Omega^+}(v)\,dv\\
\stackrel{(\star 3)}{=} & \int_{\R^2}(\rm{the\ same})\\ 
\stackrel{(\star 4)}{=} & \int_{\R^2}h_{\Omega^+}(w)\,h_{\Omega^+}(w+s-x)\,dw\\
\stackrel{(\star 5)}{=} & \int_{\Omega^+}h_{\Omega^+}(w)\,h_{\Omega^+}(w+s-x)\,dw\,,
\endaligned
$$
therefore it is the same as the first integral. 

Before proceeding let us comment on the calculation above. In $(\star 1)$ we used that $h_{\Omega^+}(\xi)=h_{\Omega^-}(-\xi)$, which comes from $\Omega$ being symmetric with respect to the origin. Next, $(\star 2)$ is simply a change of variable, whereas $(\star 3)$ and $(\star 5)$ hold because $\Omega^+$ is the support of $h_{\Omega^+}$. Finally, $(\star 4)$ follows from the substitution $w=x-s+v$. 

So we proved that 
$$
\frac{1}{|\Omega|}\int_{\Omega}\sum_{D\in {\cal G}_t}H_Df(x)\,dt=\int_{\R^2}f(s)F(x-s)\,ds\,, 
$$
where 
$$
\aligned
F(x) 
& =\frac{2}{|\Omega|} \int_{\Omega^+}h_{\Omega^+}(v)\,h_{\Omega^+}(v+x)\,dv\\
& =\frac{2}{|\Omega|} \int_{\Omega^+}h_{\Omega^+}(v)\,h_{\Omega^+}(v-x)\,dv\,.
\endaligned
$$
Note also that $\frac{2}{|\Omega|}=\frac{1}{|\Omega^+|}$ and that in case of general ${\cal H}_{a,b}^\vartriangle$ the number $|\Omega|$ equals $b$.
\qed

\bigskip
Of course, we are interested in an explicit formula for $F$. We start with the simplest case: $a=0$ and $b=1$. Then $\Omega^+$ is exactly the triangle from Figure \ref{fig0}. We actually have to compute three different kernels, since we have three different types of Haar functions. Let us denote them by $F_0$, $F_+$ and $F_-$. To them correspond (in the sense of replacing $h_{\Omega^+}$ in Proposition \ref{Piazzolla}) functions $h_0$, $h_+$ and $h_-$, respectively. They are precisely those plotted in Figures \ref{fig0}, \ref{fig+} and \ref{fig-}. These observations enable us to perform the calculation.


\setlength{\unitlength}{1mm}%
\begin{figure}
\begin{center}
\begin{picture}(140,45)(-50,0)
\put(-5,0){\vector(1,0){50}}
\put(0,-5){\vector(0,1){50}}
\put(16,16){\line(0,-1){16}}
\put(16,16){\line(-1,0){16}}
\thicklines
\put(0,0){\line(1,0){32}}
\put(0,0){\line(0,1){32}}
\put(0,32){\line(1,-1){32}}
\put(6.5,6.5){$+$}
\put(20,4){$-$}
\put(4,19){$-$}
\put(31,-4){$1$}
\put(-3,31){$1$}


\end{picture}
\end{center}
\caption{Graph of $h_{0}$}
{\protect{\label{fig0}}}%
\end{figure}
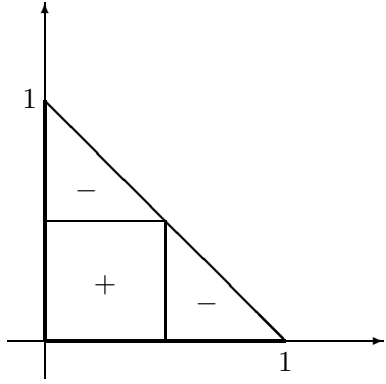



\setlength{\unitlength}{1mm}%
\begin{figure}
\begin{center}
\begin{picture}(140,45)(-50,0)
\put(-5,0){\vector(1,0){50}}
\put(0,-5){\vector(0,1){50}}
\thicklines
\put(16,0){\line(1,0){16}}
\put(0,16){\line(0,1){16}}
\put(0,32){\line(1,-1){32}}
\put(16,16){\line(0,-1){16}}
\put(16,16){\line(-1,0){16}}
\put(20,4){$-$}
\put(4,19){$+$}
\put(31,-4){$1$}
\put(-3,31){$1$}


\end{picture}
\end{center}
\caption{Graph of $h_{+}$}
{\protect{\label{fig+}}}%
\end{figure}



\setlength{\unitlength}{1mm}%
\begin{figure}
\begin{center}
\begin{picture}(140,45)(-50,0)
\put(-5,0){\vector(1,0){50}}
\put(0,-5){\vector(0,1){50}}
\put(0,16){\line(1,-1){16}}
\put(16,16){\line(0,-1){16}}
\put(16,16){\line(-1,0){16}}
\put(9,10){$-$}
\put(3,3){$+$}
\put(31,-4){$1$}
\put(-3,31){$1$}
\multiput(0,32)(7,-7){5}{\line(1,-1){4}}

\thicklines
\put(0,0){\line(1,0){16}}
\put(0,0){\line(0,1){16}}
\put(0,16){\line(1,0){16}}
\put(16,0){\line(0,1){16}}
\end{picture}
\end{center}
\caption{Graph of $h_{-}$}
{\protect{\label{fig-}}}%
\end{figure}


\bigskip
The graph of $F_0$ is symmetric with respect to both diagonals: to $x-y=0$ (because so is $h_0$) and to $x+y=0$ (take $\zeta-z$ as a new variable). 
Its support is the convex combination of the points $\pm(1,0)$, $\pm(0,1)$ and $\pm(1,-1)$ (see Figure \ref{F_0}). For the sake of transparency it is useful to introduce the function 
$$
G_0(z)=\int_\C h_0(\zeta)h_0(\zeta-z)\,dA(\zeta)\,,
$$
i.e. simply $G_0=F_0/2$. Note that $G_0(0)=\nor{h_0}_2^2=1$. 
Bearing in mind the symmetries, it suffices to get $G_0(x,y)$ on the areas A$-$G:
$$
\aligned
\rm A: & & & -(x-y)^2+(2x+2y-1)^2\\ 
\rm B: & & & -(x-y)^2  \\
\rm C: & & & -(1-x-y)^2  \\
\rm D: & & & \ \hskip 11pt (1-x)^2+2(x-y)(x-y-1)  \\
\rm E: & & & -(1-x)^2+(2x+2y-1)^2/2  \\   
\rm F: & & & -(1-x)^2  \\
\rm G: & & & \ \hskip 11pt (1-x)^2\,.  
\endaligned
$$
In an analogous way we derive $F_+$ and $F_-$.
Also, we have $G_\pm=\frac12F_\pm$. Their supports can be found on Figures \ref{F_+} and \ref{F_-}. 
The function $G_+$ is calculated as follows. 
$$
\aligned
\rm A: & & &  \ \hskip 11pt (2x+2y-1)^2\\ 
\rm D: & & &  \ \hskip 11pt (2x-1)^2-2y^2  \\
\rm E: & & & -(2x+2y-1)^2/2  \\   
\rm G: & & &  -2(1-x)^2\,.  
\endaligned
$$
Similarly, $G_-$ turns out to be 
$$
\aligned
\rm A: & & &   -(1-2x)(1-2y)+2(2x+2y-1)^2\\ 
\rm B: & & &   -(1-2x)(1-2y)\\
\rm D: & & &   \ \hskip 11pt (1-2x)(1-2y-4x)\,.  
\endaligned
$$
In order to calculate $G_\pm$ correctly one should recall that functions $h_\pm$ attain values $\{-2,2\}$ on their supports, while $h_0$ has $\{-\sqrt{2},\sqrt{2}\}$ as its range.

Now we define 
$$
F:=\sigma_0\,F_0+\sigma_+\,F_++\sigma_-\,F_-\,,
$$
where $\sigma_0,\sigma_+,\sigma_-\in S^1$. As before we can assume that $\sigma_0=1$. This $F$ is not to be mixed up with that of Proposition \ref{Piazzolla}, of course. Set also $G:=\frac12F$.


\setlength{\unitlength}{1mm}%
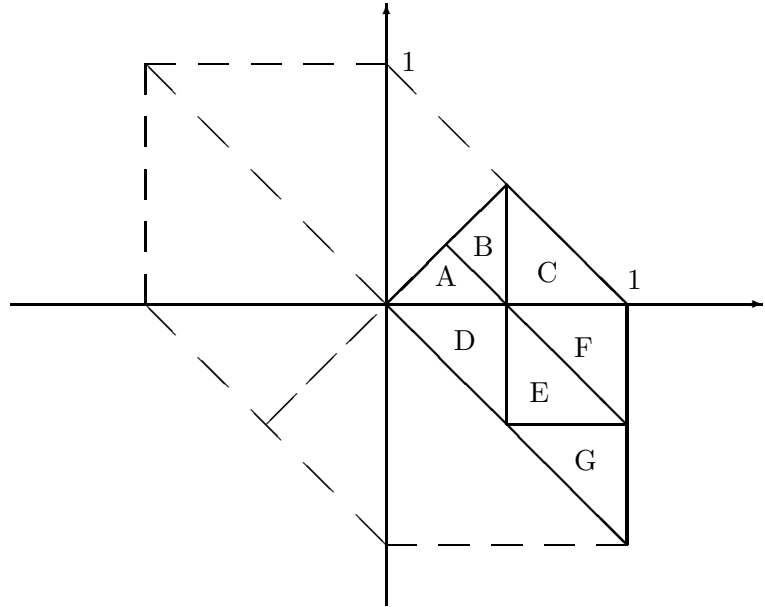
\begin{figure}
\begin{center}
\begin{picture}(140,50)(-70,-35)

\put(0,0){\vector(1,0){50}}
\put(0,0){\vector(0,1){40}}
\put(0,0){\line(-1,0){50}}
\put(0,0){\line(0,-1){40}}

\multiput(0,32)(7,-7){5}{\line(1,-1){4}}
\multiput(-32,0)(7,-7){5}{\line(1,-1){4}}
\multiput(-32,32)(7,-7){5}{\line(1,-1){4}}
\multiput(-32,32)(7,0){5}{\line(1,0){4}}
\multiput(0,-32)(7,0){5}{\line(1,0){4}}
\multiput(-32,0)(0,7){5}{\line(0,1){4}}
\multiput(-15.9,-15.9)(5.3,5.3){3}{\line(1,1){4}}



\thicklines

\put(16,16){\line(-1,-1){16}}
\put(16,16){\line(1,-1){16}}
\put(8,8){\line(1,-1){24}}
\put(16,16){\line(0,-1){32}}
\put(0,0){\line(1,-1){32}}
\put(16,16){\line(-1,-1){16}}
\put(32,0){\line(0,-1){32}}
\put(16,-16){\line(1,0){16}}

\put(32,2){$1$}
\put(2,31){$1$}
\put(6.5,2.5){A}
\put(11.5,6.5){B}
\put(20,3){C}
\put(9,-6){D}
\put(19,-13){E}
\put(25,-7){F}
\put(25,-22){G}

\end{picture}
\end{center}
\caption{Support of $G_0$}
{\protect{\label{F_0}}}%
\end{figure}



\setlength{\unitlength}{1mm}%
\begin{figure}
\begin{center}
\begin{picture}(140,50)(-70,-35)

\put(0,0){\vector(1,0){50}}
\put(0,0){\vector(0,1){40}}
\put(0,0){\line(-1,0){50}}
\put(0,0){\line(0,-1){40}}

\put(32,-0.5){\line(0,1){1}}
\put(-0.5,32){\line(1,0){1}}

\multiput(-16,32)(7,-7){5}{\line(1,-1){4}}
\multiput(-32,16)(7.35,-7.35){7}{\line(1,-1){4}}
\multiput(-32,32)(7,-7){5}{\line(1,-1){4}}
\multiput(-32,32)(6.5,0){3}{\line(1,0){3}}
\multiput(16,-32)(6.5,0){3}{\line(1,0){3}}
\multiput(-32,16)(0,6.5){3}{\line(0,1){3}}
\multiput(-7.4,-7.4)(5.3,5.3){3}{\line(1,1){4}}


\thicklines

\put(0,0){\line(1,1){8}}
\put(0,16){\line(1,-1){32}}
\put(16,0){\line(0,-1){16}}
\put(0,0){\line(1,-1){32}}
\put(32,-16){\line(0,-1){16}}
\put(16,-16){\line(1,0){16}}

\put(31,2){$1$}
\put(2,31){$1$}
\put(6.5,2.5){A}
\put(9,-6){D}
\put(19,-13){E}
\put(25,-22){G}

\end{picture}
\end{center}
\caption{Support of $G_+$}
{\protect{\label{F_+}}}%
\end{figure}
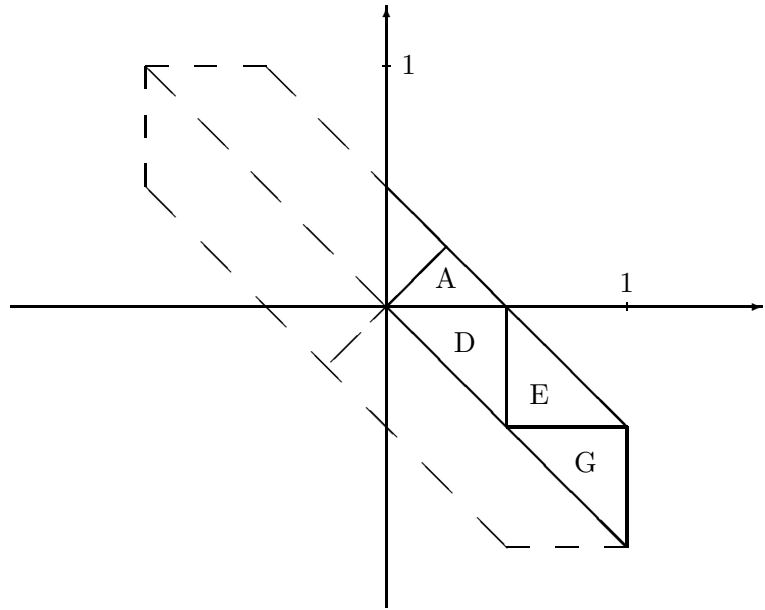



\setlength{\unitlength}{1mm}%
\begin{figure}
\begin{center}
\begin{picture}(140,50)(-70,-35)

\put(0,0){\vector(1,0){50}}
\put(0,0){\vector(0,1){40}}
\put(0,0){\line(-1,0){50}}
\put(0,0){\line(0,-1){40}}


\multiput(-16,-16)(7,0){5}{\line(1,0){3}}
\multiput(-16,16)(7.2,0){5}{\line(1,0){3}}
\multiput(-16,-16)(0,7.2){5}{\line(0,1){3}}


\put(32,-0.5){\line(0,1){1}}
\put(-0.5,32){\line(1,0){1}}

\thicklines

\put(0,0){\line(1,1){16}}
\put(16,16){\line(0,-1){32}}
\put(0,0){\line(1,-1){16}}
\put(8,8){\line(1,-1){8}}
 
\put(31,2){$1$}
\put(2,31){$1$}
\put(6.5,2){A}
\put(12,6.7){B}
\put(9,-6){D}

\end{picture}
\end{center}
\caption{Support of $G_-$}
{\protect{\label{F_-}}}%
\end{figure}
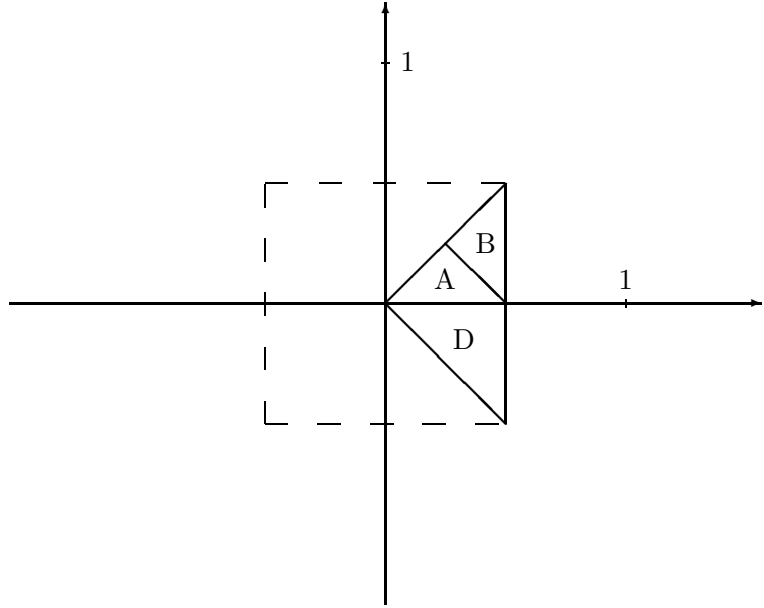


\subsubsection*{Arbitrary triangles and quest for the constant}

Next step is to bring up formulas for the kernel, corresponding to arbitrary ${\cal H}_{a,b}^\vartriangle$.
Denote by $D$ the standard triangle depicted in Figures \ref{fig0} -- \ref{fig-}. Choose $a,b\in\R$, $b\not =0$, and let $D'$ stand for the triangle on which the system ${\cal H}_{a,b}^\vartriangle$ is based. The triangle in Figure \ref{ab} represents a sufficiently general case of $D'$. Indeed, because of averaging over translations, we may assume that $D'$ has the origin as one of its verteces. Because of averaging over dilations, we may assume that one of its sides has length one. 

In order to compute the induced kernel $F'$ (the formula for which is provided by Proposition \ref{Piazzolla})
it suffices to know $F$. The map $U$, given by $U(x,y)=(x+ay,by)$, is an automorphism of $\R^2$ and also a homeomorphism $D\rightarrow D'$. Its Jacobian has determinant $b$. Also, keep in mind that the values of Haar functions will decrease by a factor of $\sqrt{b}$. Consequently 
$$
F'(x,y)=\frac1bF(U^{-1}(x,y))=\frac1bF\Big(x-\frac{a}{b}\,y,\frac{1}{b}\,y\Big)\,.
$$

Finally, we would like to find the combination of $a,b,\sigma_+,\sigma_-$ which gives the best constant after averaging. Recalling the formula \eqref{konst}, this means we have to maximize the absolute value of the integral
$$
\aligned
&\frac{1}{2\log2}\int_\R\int_\R \frac1bF\Big(x-\frac{a}{b}\,y,\frac{1}{b}\,y\big)\left(\frac{x^2-y^2}{x^2+y^2}+i\,\frac{2xy}{x^2+y^2}\right)\,dx\,dy\\
=&\frac{1}{2\log2}\int_\R\int_\R F(s,t)\frac{s+t(a+b\,i)}{s+t(a-b\,i)}\,ds\,dt
\endaligned
$$
over all $\sigma_+, \sigma_-\in S^1$, $a,b\in\R$, $b\not =0$. We may actually take only $b>0$, because changing $b$ to $-b$ essentially means conjugating the integral. Write $\zeta=a+bi$. Since $F$ was symmetric (even) with respect to the origin, this integral equals  
$$
\frac{1}{\log2}\iint_{\{s+t\geqslant 0\}} F(s,t)\frac{s+t\zeta}{s+t\bar\zeta}\,ds\,dt\,.
$$
Use also that $F$ is symmetric with respect to the diagonal $s=t$, therefore we can write the integral above as
$$
\frac{1}{\log2}\iint_{{\cal A}} F(s,t)\bigg(\frac{s+t\zeta}{s+t\bar\zeta}+\frac{t+s\zeta}{t+s\bar\zeta}\bigg)\,ds\,dt\,,
$$
where ${\cal A}=$\,A\,$\cup$\,B\,$\cup\hdots\cup$\,G (see Figure \ref{F_0}). 
Knowing that $F=2G$, the best constant will then satisfy
$$
\aligned
\frac{1}{|C|}
& =\frac{2}{\log2}
\sup_{\sigma_+, \sigma_-,a,b}\bigg|\iint_{{\cal A}} G(s,t)\bigg(\frac{s+t\zeta}{s+t\bar\zeta}+\frac{t+s\zeta}{t+s\bar\zeta}\bigg)\,ds\,dt\bigg|\,.
\endaligned
$$
The integral above is a complex number whose real and imaginary part are actually sums of up to 140 integrals of functions $G_0$, $G_+$, $G_-$ over domains A$-$G and their reflections. At this point we resort to {\sf Mathematica}. 
However, this time our numerical tests were left inconclusive, i.e. we could not arrive at any result whatsoever with a sufficient degree of certainty.

\subsubsection*{Why triangles}

All Haar functions we have been considering so far, i.e. those belonging to parallelograms, are odd or even (depending on the type of the function) with respect to the center of the parallelogram and with respect to the lines, parallel to its sides and passing through the origin. These symmetries reflect in the fact that the kernels $F$ we get are often even, whereupon the odd part of the weight ($\sin 2\varphi$) plays no r\^ole. Even when this is not the case (see Section \ref{hvd}), the best constant is attained for coefficients $\sigma$ from $\{-1,1\}$. In other words, the complex kernel ($z^{-2}$) we target is best approached with real combinations of real kernels, whoch does not seem to be very natural. 

On the other hand, triangles and their associated Haar functions do not admit such symmetries. Therefore we really cannot discard complex coefficients neither the imaginary parts of the weight $e^{2i\varphi}$. This gives rise to the hope that they might come somewhat closer to $z^{-2}*$.


We feel it might be rewarding to try in this direction. Since the group of symmetries of a triangle is smaller than the one for rectangles, we think that the ``price we pay" after averaging might be smaller, as well. In other words, this might also reflect in a constant, smaller than the one we have now.

\bigskip 

\begin{tabular}{l}
Oliver Dragi\v{c}evi\'c \\
Scuola Normale Superiore  \\
Piazza dei Cavalieri 7   \\
56126 Pisa  \\
Italy  \\ 
\\ 
{\it Current address:}   \\
Institute of Mathematics, Physics and Mechanics  \\
University of Ljubljana \\
Jadranska 19 \\
SI-1000 Ljubljana\\
Slovenia \\
oliver.dragicevic@fmf.uni-lj.si   \\
\\ \\ 
Stefanie Petermichl\\
University of Texas at Austin \\
Department of Mathematics \\
1 University Station C1200\\
Austin, TX 78712 \\
USA\\
stefanie@math.utexas.edu\\
\\ \\
Alexander Volberg   \\
Department of Mathematics  \\
Michigan State University \\
East Lansing, MI 48824 \\
USA \\ 
volberg@math.msu.edu
\end{tabular}

\newpage
\tableofcontents

\end{document}